\documentclass[3p, a4paper, 11pt]{elsarticle}

\usepackage{xcolor}
\usepackage{tikz}
\usepackage{layouts}
\biboptions{numbers,sort&compress}
\usepackage{algorithm}
\usepackage{algpseudocode}
\usepackage{gensymb}
\usepackage[unicode=true,colorlinks=true,pdfborder={0 0 0}]{hyperref}
\hypersetup{
    bookmarks=true,         
    bookmarksopen=true,     
    pdftoolbar=true,        
    pdfmenubar=true,        
    pdffitwindow=true,      
    pdfstartview={FitH},    
    pdfnewwindow=true,      
}
\usepackage{graphicx}
\usepackage{caption}
\usepackage{subcaption}
\usepackage{wrapfig}
\usepackage{multirow,makecell}
\usepackage{longtable}
\usepackage{booktabs}
\usepackage{tabularx}
\usepackage{setspace}
\usepackage{float}
\usepackage{enumerate}
\usepackage{enumitem}

\renewcommand{\appendixname}{}
\usepackage[amsthm,thmmarks]{ntheorem}
\usepackage{siunitx}
\usepackage{lineno}
\usepackage{math-symbols}
\usepackage{cases}
\usepackage{amsmath}

\graphicspath{{./images/}}

\theoremstyle{definition}
\newtheorem{theorem}{Theorem}
\newtheorem{proposition}[theorem]{Proposition}
\theoremstyle{definition}
\newtheorem{example}{Example}
\newtheorem{remark}{Remark}
\theoremstyle{definition}

\theoremstyle{definition}

\journal{X} 

\begin{document}
\begin{frontmatter}
    \title{Data-driven balanced truncation for second-order systems via the approximate Gramians}
    \author[home]{Xiaolong Wang\corref{cor1}}
    \ead{xlwang@nwpu.edu.cn}
    \author[home]{Xuerong Yang}
    \author[coop]{Xiaoli Wang}
    \author[home]{Bo Song}
    \cortext[cor1]{Corresponding author}
    \address[home]{School of Mathematics and Statistics,
        Northwestern Polytechnical University, Xi'an 710072, China}
    \address[coop]{Xi’an Microelectronics
    	Technology Institute, Xi’an 710065, China}
    \begin{abstract}
        This paper studies the data-driven balanced truncation (BT) method for 
        second-order systems based on the measurements in the frequency domain. The 
        basic idea is to approximate Gramians used the numerical quadrature rules, 
        and establish the relationship between the main quantities in the procedure 
        of BT with the sample data, which paves the way for the execution of BT in a 
        nonintrusive manner. We construct the structure-preserving reduced models 
        approximately based on the samples of second-order systems with proportional 
        damping, and provide the detailed execution of the data-driven counterpart 
        of BT in real-value arithmetic. The low-rank approximation to the solution 
        of Sylvester equations is also introduced to speed up the process of the 
        proposed approach when a large amount of samples involved in the modeling. 
        The performance of our approach is illustrated in detail via two numerical 
        examples.               
    \end{abstract}
\end{frontmatter}


\section{Introduction}\label{sec:introduction}
We consider the problem of model order reduction (MOR) for linear time-invariant 
second-order systems, which are described by the following differential equations 
\begin{equation}\label{sosystem}
	\left\{
	\begin{array}{r}
		M\ddot{q}(t)+D\dot{q}(t)+Kq(t)=Bu(t),	\\
		y(t) = Cq(t),
	\end{array}\right.
\end{equation}
where $q(t) \in \mathbb{R}^n$ is the state variable, $u(t) \in\mathbb{R}^m$ is the 
input, and $y(t) \in \mathbb{R}^p$ is the output or measurement. The coefficients  
$M$, $D$, $K\in \mathbb{R}^{n \times n}$ are the mass matrix, the damping matrix 
and the stiffness matrix, respectively. We assume that $M$ is nonsingular,
$B\in \mathbb{R}^{n\times m}$ and $C\in \mathbb{R}^{p\times n}$. With the zero 
initial conditions, 
the transfer function of (\ref{sosystem}) is defined as 
\begin{equation}\label{eq:so-transfer}
	H(s) = CG(s)B,
\end{equation}
where $G(s) = (s^2M+sD+K)^{-1}$. The focus is to construct the following 
structure-preserving reduced models of order 
$r\ll n$, such that the behavior of (\ref{sosystem}) is approximated faithfully for 
all admissible inputs
\begin{equation}\label{mor-sosystem}
	\left\{
	\begin{array}{r}
		M_r\ddot{q}_r(t)+D_r\dot{q}_r(t)+K_rq_r(t)=B_ru(t),	\\
		y_r(t) = C_rq_r(t),
	\end{array}\right.
\end{equation}
where $q_r(t) \in \mathbb{R}^r$, $M_r$, $D_r$, $K_r\in \mathbb{R}^{r\times r}$, $B_r\in 
\mathbb{R}^{r\times m}$, $C_r \in \mathbb{R}^{p\times r}$.

MOR of second-order systems has been studied intensively via a wide variety of 
approaches. In 
\cite{Bai2005,Li2012}, a second order Arnoldi method is used to generate reduced 
models with the second-order structure, instead of the linearization procedure. Some 
researchers expand second-order systems over the orthogonal polynomial basis and 
define structure-preserving reduced models via the projection matrices coming from 
the expansion coefficients \cite{Xiao2014, Eid2008}. 
The standard balanced truncation (BT) method has also been 
extended to 
the second-order case by employing or mimicking Moore's balance and truncate, and 
the detailed execution can be found in 
\cite{Meyer1996,Chahlaoui2006,reis2008,Benner2012ifac}. An alternative 
generalization of BT to second-order problems adopts the framework of generalized
Hamiltonian systems and performs the truncation by imposing a holonomic constraint 
on the system rather than standard Galerkin projection \cite{Hartmann2010}. The 
$H_2$ optimal reduced order modeling of second-order systems is reformulated as an 
optimization problem on the product manifold, and a Riemannian steepest descent 
method is exploited to generate reduced models iteratively in \cite{Sato2017}.  
The associated bitangential Hermite interpolation conditions on $H_2$ 
optimal reduction is provided in  
\cite{Mlinaric2025}. Recently, the second-order systems containing the random 
parameters are 
considered by employing a structure-preserving MOR in the framework of stochastic 
Galerkin method \cite{Pulch2024}. Basically, the second-order systems are a special 
kind of structured systems. We refer the reader to \cite{Beattie2009,Goyal2024} 
for the comprehensive coverage of MOR for structured systems in more general 
setting.

Data-driven modeling is an important tool for the simplification of complex systems. 
Several methods have been developed in the field of MOR to learn reduced models from 
data. Such kind of methods include dynamic mode decomposition 
\cite{Tu2014}, the 
Loewner framework \cite{Mayo2007,Moreschini2024}, operator inference 
\cite{Peherstorfer2016}, as well as rational 
least-squares methods AAA \cite{Rodriguez2023}. A few data-driven approaches have 
been extended for second-order systems to reserve structural features. The parameter 
optimization method is exploited for second-order systems, which calculates the 
elements of the system matrices iteratively such that the difference between the 
orginal systems and the reduced models is as small as possible 
\cite{Schwerdtner2023}. In \cite{Beattie2022}, the Loewner framework is executed to 
preserve the second order structure with internal Rayleigh-damping. The structured 
barycentric forms are proposed to model the second-order systems using the frequency 
domain input-output data in \cite{Gosea2024}, and Loewner-like algorithms are 
developed for the explicit computation of 
simplified second-order systems.

Recently, a novel reformulation of the classical BT is proposed in \cite{gosea2022}. 
It relies completely on the system response data and avoids the intrusive access to 
any prescribed realization of the original model. In this paper, we apply this basic 
idea to second-order systems, and present the nonintrusive version of BT by using 
the sample data in the frequency domain. The Gramians of second-order systems are 
first approximated via the numerical quadrature in the frequency domain, and thereby 
the main quantities involved in the standard BT can be associated with the 
measurements in the quadrature. We establish the linear equations that the main 
quantities satisfy, so as to calculate them directly via the measurements. However, 
the situation is more complex in second-order cases. There are much more variables 
in the derived linear equations, and the main quantities can not be determined 
completely. We then switch to a kind of second-order systems with proportional 
damping, where the standard BT can be reformulated as the explicit expression of the 
measurement in an approximate manner. By using the quadrature nodes and weights in a 
symmetric manner, we present a real-valued algorithm for the proposed data-driven 
BT, leading to real-valued reduced models as well. When the data-driven BT is 
implemented based on a large amount of measurements, the SVD involved in the 
algorithm is extremely time consuming. We provide low-rank approximate solutions to 
the Sylvester equations by using the extended Krylov subspace methods. Consequently, 
we just need to perform the SVD of a low order matrix in the procedure of modeling, 
thereby enabling an efficient execution of our approach.

The paper is organized as follows. \autoref{sec:sec-2} introduces the standard BT 
and the preliminaries on Gramians. We start \autoref{sec:sec-3} with 
the approximation to Gramians, and then establish the data-driven BT for 
second-order systems with proportional damping in real arithmetic. A low-rank 
approximation to the solutions of Sylvester equations is also provided to enable an 
efficient execution of our approach. Numerical examples are 
used to test our approach in \autoref{sec:sec-4}.
Finally, some conclusions are drawn in \autoref{sec:sec-5}.

\section{Preliminaries on BT for second-order systems}\label{sec:sec-2}

In this section, we briefly review BT methods for second-order systems to facilitate 
the description of our data-driven approach. Gramians of second-order systems are 
defined via the equivalent first-order systems. With 
the new state 
$x(t)=\left[\begin{array}{cc}q^\top(t) & \dot{q}^\top(t)\end{array}\right]^\top$, 
second-order systems (\ref{sosystem}) can be reformulated as the following linear system 
\begin{equation}\label{fosystem}
	\left\{
	\begin{array}{l}
		\mathcal{E} \dot{x}(t)=\mathcal{A} x(t)+\mathcal{B} u(t),	\\
		y(t)=\mathcal{C} x(t),
	\end{array}\right.
\end{equation}
where the coefficient matrices are 
\begin{equation*}
	\begin{array}{cccc}
	\mathcal{E} = \left[\begin{array}{rr}
		I &  0\\
		0 &  M
	\end{array}\right],	& \mathcal{A} = \left[\begin{array}{rr}
	0  & I \\
	-K & -D 
	\end{array}\right], & \mathcal{B} = \left[\begin{array}{r}
	0\\
	B
	\end{array}\right], & \mathcal{C}=\left[\begin{array}{rr}
	C& 0 
	\end{array}\right].
	\end{array}
\end{equation*}
By the standard BT approach, the controllability and observability Gramians of 
(\ref{fosystem}) are defined explicitly as follows
\begin{equation*}
	\begin{split}
	\mathcal{P}=\int_{0}^{\infty} e^{\mathcal{E}^{-1} \mathcal{A} t} \mathcal{E}^{-1} 
	\mathcal{B B}^{\top} \mathcal{E}^{-\top} e^{\mathcal{A}^{\top} \mathcal{E}^{-\top} t} 
	dt,\\
	\mathcal{Q}=\int_{0}^{\infty} \mathcal{E}^{-\top} e^{\mathcal{A}^{\top} 
	\mathcal{E}^{-\top} t} \mathcal{C}^\top\mathcal{C}e^{\mathcal{E}^{-1} \mathcal{A} t} 
	\mathcal{E}^{-1} dt.
	\end{split}
\end{equation*}
By Parseval's theorem, $\mathcal{P}$ and $\mathcal{Q}$ can be represented in the frequency domain as
\begin{equation}\label{frP}
	\mathcal{P}=\frac{1}{2 \pi} \int_{-\infty}^{\infty}(i \xi 
	\mathcal{E}-\mathcal{A})^{-1} \mathcal{B B}^{\top}\left(-i \xi  
	\mathcal{E}^{\top}-\mathcal{A}^{\top}\right)^{-1} d \xi,
\end{equation}
\begin{equation}\label{frQ}
	\mathcal{Q}=\frac{1}{2 \pi} \int_{-\infty}^{\infty}\left(-i \omega 
	\mathcal{E}^{\top}-\mathcal{A}^{\top}\right)^{-1} \mathcal{C}^{\top} \mathcal{C}(i 
	\omega \mathcal{E}-\mathcal{A})^{-1} d \omega.
\end{equation}
Let Gramians (\ref{frP}), (\ref{frQ}) be partitioned as
\begin{equation}
 	\begin{array}{ccc}
 	\mathcal{P} = \left[\begin{array}{cc}
 	 P_p	&  P_{pv} \\
 	 P^\top_{pv}	& P_v
 	\end{array}\right]	&\text{and} &\mathcal{Q} = \left[\begin{array}{cc}
 	Q_p     & Q_{pv} \\
 	Q^\top_{pv}    & Q_v
 \end{array}\right],
 	\end{array}
\end{equation}
with all the blocks being of size $ n \times n$. Then $P_p$ and $P_v$ are defined as the 
position and velocity controllability Gramians, while $Q_p$ and $Q_v$ are the position 
and 
velocity observability Gramians of second-order systems (\ref{sosystem}). By the 
block structure of matrices in (\ref{frP}) and (\ref{frQ}), we perform the basic matrix 
manipulation and obtain the following 
explicit expression of Gramians 
\begin{eqnarray}
	&&P_p = \frac{1}{2\pi}\int_{-\infty}^{\infty}G(i\xi)BB^\top G^{\mathrm H}(i\xi)d\xi, 
	\label{G-1}
	\\
	&&P_v = \frac{1}{2\pi}\int_{-\infty}^{\infty}G(i\xi)\xi^2BB^\top G^{\mathrm 
	H}(i\xi)d\xi, \label{G-2}\\
	&&Q_p = \frac{1}{2\pi}\int_{-\infty}^{\infty}(-i\omega M^\top+D^\top)G^{\mathrm 
	H}(i\omega)C^\top CG(i\omega)(i\omega M+D)d \omega,\label{G-3}\\
	&&Q_v = \frac{1}{2\pi}\int_{-\infty}^{\infty}G^{\mathrm H}(i\omega)C^\top 
	CG(i\omega)d \omega,\label{G-4}
\end{eqnarray}
where $G(\cdot)$ is defined as in transfer function (\ref{eq:so-transfer}). In practice, 
the numerical approximation to Gramians (\ref{G-1})-(\ref{G-4}) is 
typically obtained by solving the associated Lyapunov equations satisfied by $\mathcal P$ 
and $\mathcal Q$. We refer the reader to \cite{Chahlaoui2006} for more details.

The procedure of BT for second-order systems can be conducted based on the different 
balanced realizations, corresponding to the individual singular values of 
(\ref{sosystem}). For example, the velocity singular values are defined as the square 
roots of the eigenvalues of the matrix $P_vM^\top Q_vM$, while the position-velocity 
singular values are the square roots of the eigenvalues of the matrix $P_pM^\top Q_vM$. 
The balanced transformation of (\ref{sosystem}) is designed to make the controllability 
and observability Gramians equal and diagonal 
\begin{equation*}
	\begin{split}
		\hat{P}_v = \hat{Q}_v &= \Sigma_v\quad \text{(velocity balanced)}, \\
		\hat{P}_p = \hat{Q}_v &= \Sigma_{pv}\quad \text{(position-velocity balanced)}. 
	\end{split}
\end{equation*}
The singular values allow us to determine the important positions and velocities, i.e., 
those with large effect on the input-output map. Reduced models (\ref{mor-sosystem}) can 
be derived by the 
transformation of (\ref{sosystem}) into one kind of balanced forms, followed by the 
direct truncation of the 
positions and velocities corresponding to the small singular values.
In this paper, we take the velocity balanced truncation as an example to elaborate 
the data-driven BT of second-order systems. Based on the square factors $L$, $U \in 
\mathbb{R}^{n \times n}$ such that 
\begin{equation*}
	\begin{array}{cc}
	P_v = UU^T,  & Q_v = LL^T,
	\end{array}
\end{equation*}
we summarize in Algorithm \ref{SOBT} the velocity BT procedure of second-order systems 
\cite{reis2008}. Based on the observation that the main terms in step 2 and 
step 4 of Algorithm \ref{SOBT} may be approximated by the samples of the transfer 
function, a data-driven approach will be presented in next section.  

\begin{algorithm}[htb]
	\caption{Velocity balanced and truncated method for second order systems.}
	\label{SOBT}
	\begin{algorithmic}[1]
		\Require
		System matrices $M$, $D$, $K$, $B$, $C$;
		\Ensure
		Reduced models $M_r\in\mathbb{R}^{r \times r}$, $D_r\in\mathbb{R}^{r\times r}$, 
		$K_r\in\mathbb{R}^{r\times r}$, $B_r\in\mathbb{R}^{r\times m}$, 
		$C_r\in\mathbb{R}^{p\times r}$;
		\State Compute the square factors $P_v = UU^{\top}$, $Q_v = LL^{\top}$, and 
		pick a 
		truncation index, $1\leq r \leq \min(\mathrm{rank}(U), \mathrm{rank}(L))$.
		\State Compute the SVD of the matrix $L^{\top}MU$, with the partitioned form 
		as 
		follows
		\begin{equation*}
			L^{\top}MU=\left[\begin{array}{cc}
				Z_1	& Z_2
			\end{array}\right]\left[\begin{array}{cc}
				S_1&  \\
				& S_2
			\end{array}\right]\left[\begin{array}{c}
				Y^{\top}_1\\
				Y^{\top}_2
			\end{array}\right],
		\end{equation*}
		where $S_1\in \mathbb{R}^{r\times r}$ and $S_2\in \mathbb{R}^{(n-r)\times(n-r)}$.
		\State Construct the reduction basis matrices 
		\begin{equation*}
			\begin{array}{cc}
				W_r=LZ_1S^{-1/2}_1	&\text{and}\quad V_r=UY_1S^{-1/2}_1
			\end{array}.
		\end{equation*}
		\State The reduced models are given by
		\begin{equation*}
			\begin{array}{c}
				M_r=W^{\top}_rMV_r=I_r,	\\
				D_r=W^{\top}_rDV_r=S_{1}^{-1 / 2}Z^{\top}_1(L^{\top}DU)Y_1S_{1}^{-1 
				/ 
				2},	\\
				K_r=W^{\top}_rKV_r=S_{1}^{-1 / 2}Z^{\top}_1(L^{\top}KU)Y_1S_{1}^{-1 
				/ 
				2},	\\
				B_r=W^{\top}_rB=S_{1}^{-1 / 2}Z^{\top}_1(L^{\top}B),	\\
				C_r=CV_r=(CU)Y_1S_{1}^{-1 / 2}.
			\end{array}
		\end{equation*}
	\end{algorithmic}
\end{algorithm}

\section{Data-driven BT for second-order systems}\label{sec:sec-3}

We present a data-driven BT framework for second-order systems. Unlike the standard BT 
approach, it relies on the 
sample data in the frequency domain, and does not require to access to the state-space 
realization of original systems.
For the purpose of brevity, we first focus on the single-input 
single-output (SISO) case, i.e., 
$B \in \mathbb{R}^{n \times 1},C \in \mathbb{R}^{1\times n}$, while the results can be 
extended to the multiple-input multiple-output (MIMO) case as well with some proper 
modification.

\subsection{Approximation to the main quantities in BT procedure}
We consider a numerical quadrature rule to approximate the velocity controllability 
Gramian matrix
\begin{equation}\label{P-app}
	P_v \approx \widetilde{P}_v=\sum_{j=1}^{N_p}\rho^2_jG(i\zeta_j)\zeta^2_jBB^\top 
	G^\mathrm{H}(i\zeta_j),
\end{equation}
where 
$\rho^2_j$ and $\zeta_j$
represent the numerical quadrature weights and nodes, respectively, and $N_p$
is the total number of quadrature nodes. Based on this expression, the Gramian matrix 
$\widetilde{P}_{v}$
can be decomposed as $\widetilde{P}_{v}=\widetilde{U}\widetilde{U}^\mathrm{H}$
, where 
$\widetilde{U}$
is the square-root factor of the matrix
\begin{equation}\label{U}
	\widetilde{U}=\left[\begin{array}{ccc}
		\rho_1\zeta_1G(i\zeta_1)B, &\dots, &\rho_{N_p}\zeta_{N_p}G(i\zeta_{N_p})B
	\end{array}\right]\in \mathbb{C}^{n \times N_p}.
\end{equation}
Similarly, the velocity observability Gramian matrix $Q_v$ can also be approximated by 
the following numerical quadrature
\begin{equation}\label{Q-app}
	Q_v \approx \widetilde{Q}_v=\sum_{k=1}^{N_q}\varphi^2_kG^\mathrm{H}(i\omega_k)C^\mathrm{T}CG(i\omega_k),
\end{equation}
where $\varphi^2_k$ and $\omega_k$ represent the numerical quadrature weights and nodes, 
respectively, and $N_q$ is the total number of quadrature nodes. Let $\widetilde{Q}_v = 
\widetilde{L}\widetilde{L}^\mathrm{H}$. The corresponding square-root factor is given as
\begin{equation}\label{L}
	\widetilde{L}^\mathrm{H}=\left[\begin{array}{c}
		\varphi_1CG(i\omega_1)	\\
		\vdots\\
		\varphi_{N_q}CG(i\omega_{N_q})
	\end{array}\right]\in \mathbb{C}^{N_q \times n}.
\end{equation}

For the convenience of presentation,
we use the new notations for the following quantities 
\begin{equation*}
\widetilde{\mathbb{M}} = \widetilde{L}^\mathrm{H}M\widetilde{U}, \widetilde{\mathbb{D}} 
= \widetilde{L}^\mathrm{H}D\widetilde{U}, \widetilde{\mathbb{K}} = 
\widetilde{L}^\mathrm{H}K\widetilde{U}, \widetilde{\mathbb{B}} = 
\widetilde{L}^\mathrm{H}B, \widetilde{\mathbb{C}} = C\widetilde{U}. 
\end{equation*}
Clearly, there holds
\begin{equation*}
	\widetilde{\mathbb{B}} = \left[\begin{array}{c}
		\varphi_1H(i\omega_1)	\\
		\vdots\\
		\varphi_{N_q}H(i\omega_{N_q})
	\end{array}\right] \in \mathbb{C}^{N_q\times 1},
	\widetilde{\mathbb{C}} = \left[\begin{array}{ccc}
		\rho_1\zeta_1H(i\zeta_{1}),	& \dots, & \rho_{N_p}\zeta_{N_p}H(i\zeta_{N_p})
	\end{array}\right]\in \mathbb{C}^{1 \times N_p},
\end{equation*}
where the elements are 
\begin{equation}\label{LBCU}
		 \widetilde{\mathbb{B}}_k = \varphi_k H(i\omega_k),  
		\widetilde{\mathbb{C}}_j = \rho_j \zeta_j H(i\zeta_j), \quad k = 1, \cdots, N_q, 
		j = 1, \cdots, N_p. 
\end{equation}
That is, the matrices $\widetilde{\mathbb{B}}$ and $\widetilde{\mathbb{C}}$ can be 
derived 
directly via the samples of $H(s)$ for a given numerical quadrature rule. For the new 
defined matrices 
$\widetilde{\mathbb{M}}$, $\widetilde{\mathbb{D}}$, $\widetilde{\mathbb{K}}\in 
\mathbb{C}^{N_q \times N_p}$, we have the following proposition.

\begin{proposition}\label{pro:3}
	Let $\widetilde{U}$
	and $\widetilde{L}$ be the matrices defined in \textup{(\ref{U})} and 
	\textup{(\ref{L})}, respectively. Then the matrices 
	$\widetilde{\mathbb{M}}$, $\widetilde{\mathbb{D}}$, $\widetilde{\mathbb{K}}\in 
	\mathbb{C}^{N_q \times N_p}$ satisfy the following linear equation 	
	\begin{equation}\label{eqMDK1}
		\Lambda_{Q}^{2} \widetilde{\mathbb{M}}+\Lambda_{Q} 
		\widetilde{\mathbb{D}}+\widetilde{\mathbb{K}}=\Phi_{Q}^{\mathrm{T}} 
		\overrightarrow{H}(i \zeta),
	\end{equation}
	\begin{equation}\label{eqMDK2}
		\widetilde{\mathbb{M}} \Lambda_{P}^{2}+
		\widetilde{\mathbb{D}}\Lambda_{P}+\widetilde{\mathbb{K}}=\overrightarrow{H}(i 
		\omega)^{\mathrm{T}} \Xi_{P},
	\end{equation}
	where the coefficient matrices are defined as 
	\begin{equation*}
		\Lambda_{Q} = \mathrm{diag}\{i\omega_1, \cdots, i\omega_{N_q}\}\in\mathbb{C}^{N_q 
		\times N_q}, \Lambda_{P} = \mathrm{diag}\{i\zeta_1, \cdots, 
		i\zeta_{N_p}\}\in\mathbb{C}^{N_p \times N_p},
	\end{equation*}
	\begin{equation*}
		\Phi_{Q} = \left[\begin{array}{ccc}
			\varphi_1 & \cdots & \varphi_{N_q} 
		\end{array}\right] \in \mathbb{C}^{1 \times N_q},\quad 
		\Xi_{P} = \left[\begin{array}{ccc}
			\rho_1\zeta_1 & \cdots & \rho_{N_p}\zeta_{N_p}
		\end{array}\right]\in \mathbb{C}^{1 \times N_p},
	\end{equation*}
	\begin{equation*}
		\overrightarrow{H}(i \zeta) = \left[\begin{array}{ccc}
			\rho_1\zeta_1H(i\zeta_1) &  \cdots  & \rho_{N_p}\zeta_{N_p}H(i\zeta_{N_p})
		\end{array}\right]\in \mathbb{C}^{1 \times N_p},
	\end{equation*}
	\begin{equation*}
		\overrightarrow{H}(i \omega) = \left[\begin{array}{ccc}
			\varphi_1H(i\omega_1) & \cdots  & \varphi_{N_q}H(i\omega_{N_q})
		\end{array}\right]\in \mathbb{C}^{1\times N_q}.
	\end{equation*}

	\begin{proof}
		The direct matrix manipulation leads to
		\begin{equation*}
			\begin{split}
			\widetilde{\mathbb{M}} &= \widetilde{L}^\mathrm{H}M\widetilde{U}\\ 
			&= 
			\left[\begin{array}{c}
				\varphi_1CG(i\omega_1)	\\
				\vdots\\
				\varphi_{N_q}CG(i\omega_{N_q})
			\end{array}\right]M\left[\begin{array}{ccc}
				\rho_1\zeta_1G(i\zeta_1)B,	 \dots , 
				\rho_{N_p}\zeta_{N_p}G(i\zeta_{N_p})B
			\end{array}\right]\\
			&=\left[\begin{array}{ccc}
				\varphi_1\rho_1\zeta_1CG(i\omega_1)MG(i\zeta_1)B	& \dots & 
				\varphi_1\rho_{N_p}\zeta_{N_p}CG(i\omega_1)MG(i\zeta_{N_p})B \\
				\vdots	& \ddots &  \vdots\\
				\varphi_{N_q}\rho_1\zeta_1CG(i\omega_{N_q})MG(i\zeta_1)B	& \dots & 
				\varphi_{N_q}\rho_{N_p}\zeta_{N_p}CG(i\omega_{N_q})MG(i\zeta_{N_p})B
			\end{array}\right]. 
		\end{split}
		\end{equation*}
The $(k,j)$ entry of $\widetilde{\mathbb M}$ reads
$$\widetilde{\mathbb{M}}_{k,j}=\varphi_{k}\rho_{j}\zeta_{j}CG(i\omega_k)MG(i\zeta_j)B, 
\quad 1 \leq k \leq N_q, 1 \leq j \leq N_p. $$ Similarly, for 
$\widetilde{\mathbb D}$ and $\widetilde{\mathbb K}$, we have 
$$\widetilde{\mathbb{D}}_{k,j}=\varphi_{k}\rho_{j}\zeta_{j}CG(i\omega_k)DG(i\zeta_j)B 
\quad\text{and} \quad
\widetilde{\mathbb{K}}_{k,j}=\varphi_{k}\rho_{j}\zeta_{j}CG(i\omega_k)KG(i\zeta_j)B.$$
One can verify the identity 
\begin{equation*}
	\begin{split}
	&\varphi_{k}\rho_{j}\zeta_{j}H(i\zeta_j)-\varphi_{k}\rho_{j}\zeta_{j}H(i\omega_k)\\
	=&\varphi_{k}\rho_{j}\zeta_{j}CG(i\omega_k)\left[(-\omega^2_kM+i\omega_kD+K)-
	(-\zeta^2_jM+i\zeta_jD+K)\right]G(i\zeta_j)B\\
	=&(-\omega^2_k+\zeta^2_j)\widetilde{\mathbb{M}}_{k,j}+(i\omega_k-i\zeta_j)\widetilde{\mathbb{D}}_{k,j}.
	\end{split}
\end{equation*}
Similarly, there holds another identity 
\begin{equation*}
	\begin{split}
&\varphi_{k}\rho_{j}\zeta_{j}H(i\zeta_j)+\varphi_{k}\rho_{j}\zeta_{j}H(i\omega_k)\\
=&\varphi_{k}\rho_{j}\zeta_{j}CG(i\omega_k)\left[(-\omega^2_kM+i\omega_kD+K)+(-\zeta^2_jM+i\zeta_jD+K)\right]G(i\zeta_j)B\\
=&(-\omega^2_k-\zeta^2_j)\widetilde{\mathbb{M}}_{kj}+(i\omega_k+i\zeta_j)
\widetilde{\mathbb{D}}_{k,j}+2\widetilde{\mathbb{K}}_{k,j}.
	\end{split}
\end{equation*}
The above identities immediately lead to 
		\begin{equation*}
			\begin{split}
				-\omega^2_k\widetilde{\mathbb{M}}_{k,j} + 
				i\omega_k\widetilde{\mathbb{D}}_{k,j} + \widetilde{\mathbb{K}}_{k,j}&= 
				\varphi_{k}\rho_{j}\zeta_{j}H(i\zeta_j), \\
				-\zeta^2_j\widetilde{\mathbb{M}}_{k,j} + 
				i\zeta_j\widetilde{\mathbb{D}}_{k,j} + \widetilde{\mathbb{K}}_{k,j}&= 
				\varphi_{k}\rho_{j}\zeta_{j}H(i\omega_k),
			\end{split}
		\end{equation*}
which are exactly the $(k,j)$ element of (\ref{eqMDK1}) and (\ref{eqMDK2}), respectively. 
It concludes the proof. 
	\end{proof}
\end{proposition}

\begin{remark}
 The approximation to Gramians relies on the specific quadrature rule. For (\ref{P-app}) 
 and (\ref{Q-app}), there is no “node at infinity”  involved in the summation for the 
 ease of presentation. However, the potential “node at infinity” may occur if the 
 Clenshaw-Curtis quadrature is employed in the approximation \cite{Boyd1987}. 
\end{remark}

\begin{remark}
As the velocity singular values are the square roots of the eigenvalues of 
$P_vM^\top Q_vM$, they are exactly the singular values of the matrix $L^\top MU$ 
because of the square-root decomposition $P_v=UU^\top$ and $Q_v=LL^\top$. The error 
induced by the numerical quadrature can be bounded. Let 
$\sigma_1\geq\sigma_2\geq\cdots\geq\sigma_n$ denote singular values of $L^TMU$, and 
$\widetilde{\sigma}_1 \geq \widetilde{\sigma}_2 \geq \dots \geq 
\widetilde{\sigma}_{n}$ denote the singular values of 
$\widetilde{L}^\mathrm{H}M\widetilde{U}$. If 
$$\parallel 
Q_{v}-\widetilde{Q}_{v}\parallel_{F}\leq\frac{\delta}{1+\delta}\sigma_{min}(Q_{v}), 
\parallel 
P_{v}-\widetilde{P}_{v}\parallel_{F}\leq\frac{\delta}{1+\delta}\sigma_{min}(P_{v})$$
for some $\delta \in (0,1)$, where $\sigma_{min}$ denotes the smallest singular 
value, then there holds
\begin{equation*}
	\left(\sum_{k=1}^n(\sigma_k-\widetilde{\sigma}_k)^2\right)^{\frac{1}{2}}\le2\delta\parallel
	 M\parallel_2\parallel L\parallel_2\parallel U\parallel_2.
\end{equation*}
The error bound can be derived directly from the analysis on the Hankel singular values 
in \cite{gosea2022}, and we omit the details on the proof. 
\end{remark}

\begin{remark}
	Note that when the coefficient matrix $M=0$, (\ref{sosystem}) boils down to 
	a typical linear system, and linear equations (\ref{eqMDK1}) and (\ref{eqMDK2}) 
	reduce to
	the counterpart of the BT procedure for general linear systems. In 
	\cite{gosea2022}, the 
	explicit expressions for the main quantities in BT are obtained with some 
	sophisticated tricks, and one can validate that they exactly satisfy the above 
	linear equations. 
	However, the matrices 
	$\widetilde{\mathbb M},\widetilde{\mathbb D}$ and $\widetilde{\mathbb K}$ can 
	not be determined completely by (\ref{eqMDK1}) and (\ref{eqMDK2}) in the case 
	of second-order systems. We will turn to second-order systems with proportional 
	damping to achieve a data-driven scheme in next subsection. 
\end{remark}

\subsection{Data-driven BT of second-order systems with proportional damping}
We consider the proportional damping
hypothesis for second-order systems, i.e., the damping matrix $D$ is given by a linear 
combination of the mass
and stiffness matrices
\begin{equation}\label{damping-relation}
	D = \alpha M + \beta K,
\end{equation} 
for $\alpha$, $\beta \geq 0$. This specific choice of the damping matrix is often 
exploited 
in various engineering application \cite{Wu2015,Meirovitch1997}. We refer the reader 
to 
\cite{Beattie2005pro,Bonin2016,Beattie2022} for more details on MOR of this 
special case.

With the information coming from 
(\ref{damping-relation}), linear systems (\ref{eqMDK1}) and (\ref{eqMDK2}) become 
\begin{equation}\label{relation_damp}
	\begin{split}
		\Lambda_{Q}^{2} \widetilde{\mathbb{M}}+\Lambda_{Q} 
		(\alpha \widetilde{\mathbb{M}} + \beta 
		\widetilde{\mathbb{K}})+\widetilde{\mathbb{K}}&=\Phi_{Q}^{\mathrm{T}}\overrightarrow{H}(i
		 \zeta),\\		
		\widetilde{\mathbb{M}} 
		\Lambda_{P}^{2}+(\alpha \widetilde{\mathbb{M}} + \beta 
		\widetilde{\mathbb{K}})\Lambda_{P}+\widetilde{\mathbb{K}}&=\overrightarrow{H}(i
		 \omega)^{\mathrm{T}} \Xi_{P}. 
	\end{split}
\end{equation}
In order to determine the matrices $\widetilde{\mathbb{M}},\widetilde{\mathbb{K}}$ 
completely, we assume that the quadrature nodes are disjoint for the numerical quadrature 
of $P_v$ and $Q_v$, that is, $\omega_k\neq\zeta_j$ for $1 \leq k \leq N_q$ and $1 \leq j 
\leq N_p$. Specifically, the 
$(k,j)$ element of the matrices $\widetilde{\mathbb{M}}$, $\widetilde{\mathbb{D}}$, 
$\widetilde{\mathbb{K}} \in \mathbb{C}^{N_q \times N_p}$ have the following explicit 
expression 
	\begin{equation}\label{M_exp}
		\widetilde{\mathbb{M}}_{k,j}=\frac{\varphi_k\rho_j\zeta_j(H(i\zeta_j)-H(iw_k))}{\alpha(\omega_k)-\alpha(\zeta_j)}
		+\frac{\varphi_k\rho_j\zeta_j(\beta(\zeta_j)-\beta(\omega_k))(\alpha(\omega_k)H(i\omega_k)-\alpha(\zeta_j)H(i\zeta_j))}{(\alpha(\omega_k)-\alpha(\zeta_j))(\alpha(\omega_k)\beta(\zeta_j)-\alpha(\zeta_j)\beta(\omega_k))},
	\end{equation}
	\begin{equation}\label{K_exp}
		\widetilde{\mathbb{K}}_{k,j}=\frac{\varphi_k\rho_j\zeta_j\left(\alpha(\omega_k)H(i\omega_k)-\alpha(\zeta_j)H(i\zeta_j)\right)}{\alpha(\omega_k)\beta(\zeta_j)-\alpha(\zeta_j)\beta(\omega_k)},
	\end{equation}
	\begin{equation}\label{D_exp}
		\widetilde{\mathbb{D}}_{k,j}=\alpha\widetilde{\mathbb{M}}_{k,j}+\beta\widetilde{\mathbb{K}}_{k,j},
	\end{equation}
	where $\alpha(\zeta)=(-\zeta^2+i\zeta \alpha)$, $\beta(\zeta)=(i\zeta \beta+1)$.

	Now we are in a position to present the data-driven BT method for second order 
	systems. With the approximation (\ref{P-app}) and (\ref{Q-app}), the main terms 
	$L^{\top}MU, 
	L^{\top}DU,L^{\top}KU,L^{\top}B$ and $CU$ in Algorithm \ref{SOBT} can be 
	calculated approximately via 
	the samples of the frequency response, as shown in the explicit expression of 
	$\tilde{\mathbb{M}}, \tilde{\mathbb{D}}, \tilde{\mathbb{K}}, 
	\tilde{\mathbb{B}}$ and $\tilde{\mathbb{C}}$. The main steps of the proposed 
	method are summarized in 
	Algorithm \ref{al:dataSOBT}.
	\begin{algorithm}[htb]
		\caption{Data-driven BT of second-order systems with 
			proportional damping. }
		\label{al:dataSOBT}
		\begin{algorithmic}[1]
			\Require		
			Quadrature nodes $\zeta_j, \omega_k$ and weights $\rho_j, 
			\varphi_k$, for $j = 1\dots N_p, k = 1\dots 
			N_q$;
			
			Sample data of $H(s)$ at the quadrature nodes, and the 
			index $1 \leq r \leq \mathrm{min}(N_p, N_q)$.
			
			\Ensure
			Reduced models $M_r\in\mathbb{R}^{r\times 
				r},D_r\in\mathbb{R}^{r\times r},K_r\in\mathbb{R}^{r\times 
				r},B_r\in\mathbb{R}^{r\times 1},C_r\in\mathbb{R}^{1\times r}$. 
			
			\State Assemble the data  
			$\left\{H(i\zeta_j)\right\}^{N_p}_{j=1}$
			and $\left\{H(i\omega_k)\right\}^{N_q}_{k=1}$, and calculate the main 
			terms $\widetilde{\mathbb{M}}$, 
			$\widetilde{\mathbb{K}}$, $\widetilde{\mathbb{D}}$, 
			$\widetilde{\mathbb{B}}$, $\widetilde{\mathbb{C}}$ by (\ref{M_exp}), 
			(\ref{K_exp}), (\ref{D_exp}), and (\ref{LBCU}), respectively.
			
			\State Compute the SVD of the matrix $\widetilde{\mathbb{M}}$
			\begin{equation*}
				\widetilde{\mathbb{M}}=\left[\begin{array}{cc}
					\widetilde{Z}_1	& \widetilde{Z}_2
				\end{array}\right]\left[\begin{array}{cc}
					\widetilde{S}_1&  \\
					& \widetilde{S}_2
				\end{array}\right]\left[\begin{array}{c}
					\widetilde{Y}^\mathrm{H}_1\\
					\widetilde{Y}^\mathrm{H}_2
				\end{array}\right],
			\end{equation*}
			where $\widetilde{S}_1\in \mathbb{R}^{r\times r}$ and the factors are 
			partitioned compatibly.

			\State The reduced models are given by
			\begin{equation*}
				M_r=I_r, 
				D_r=\widetilde{S}^{-1/2}_1\widetilde{Z}^\mathrm{H}_1\widetilde{\mathbb{D}}
				\widetilde{Y}_1\widetilde{S}^{-1/2}_1,
				K_r=\widetilde{S}^{-1/2}_1\widetilde{Z}^\mathrm{H}_1\widetilde{\mathbb{K}}
				\widetilde{Y}_1\widetilde{S}^{-1/2}_1,
				B_r=\widetilde{S}^{-1/2}_1\widetilde{Z}^\mathrm{H}_1\widetilde{\mathbb{B}},
				C_r=\widetilde{\mathbb{C}}\widetilde{Y}_1\widetilde{S}^{-1/2}_1.
			\end{equation*}
		\end{algorithmic}
	\end{algorithm}

In practice, the dynamical systems are determined by the real-value coefficient 
matrices in general, 
ensuring that the real-valued inputs and initial conditions result in the 
real-valued outputs. It is preferred to generate a real-valued reduced model as 
well. Unfortunately, it is clear that the calculation in Algorithm \ref{al:dataSOBT} 
involves the complex arithmetic and the resulting reduced models are determined 
typically by 
complex-valued matrices. In fact, one can produce produced real-valued 
reduced models by selecting the quadrature nodes and weights in a symmetric manner 
for the numerical integration. We follow the strategy provided in \cite{gosea2022} 
with some proper modification to achieve this goal.

Let the number of quadrature nodes in both 
sets be even, that is $N_p = 2\nu_p$ and $N_q = 2\nu_q$. We assume that the 
quadrature 
nodes are symmetrically 
distributed along the real axis, i.e., 
	\begin{equation*} 
		\begin{split}
	\zeta_{-\nu_p}<\zeta_{-\nu_p+1}<\cdots<\zeta_{-1}<0<\zeta_1<\cdots<
	\zeta_{\nu_p-1}<\zeta_{\nu_p},\\
	\omega_{-\nu_q}<\omega_{-\nu_q+1}<\cdots<\omega_{-1}<0<\omega_1<\cdots
	<\omega_{\nu_q-1}<\omega_{\nu_q},
	 	\end{split}
	\end{equation*} 
such that $\zeta_j = -\zeta_{-j}, \omega_k = -\omega_{-k}$ and the corresponding 
weights $\rho_j = 
	\rho_{-j}, \varphi_k = \varphi_{-k}$ for $j = 1, \dots, \nu_p, k = 1, \dots, 
	\nu_q$. By rearranging the nodes 
	and weights, we obtain 
		\begin{equation}\label{order-nodes}
		\{\zeta_{1},\zeta_{-1},\zeta_{2},\zeta_{-2},\ldots,\zeta_{\nu_{p}},\zeta_{-\nu_{p}}\}\quad\text{and}\quad\{\omega_{1},\omega_{-1},\omega_{2},\omega_{-2},\ldots,\omega_{\nu_{q}},\omega_{-\nu_{q}}\},
		\end{equation} 
where the quadrature nodes appear in the series in pairs of conjugation. Note that 
the evaluation of the transfer function satisfies the following relationship 
\begin{equation*} 
H(\overline{s})=C^\mathsf{T}(\overline{s}^2M+\overline{s}D+K)^{-1}B=
	\overline{C^\mathsf{T}(s^2M+sD+K)^{-1}B}=\overline{H(s)},
\end{equation*} 	
where $\overline{s}$ is the conjugation of the complex number $s$. 
As a result, when the evaluation of the transfer function in Algorithm 
\ref{al:dataSOBT} is 
reordered in the same way, they form a series of conjugate pairs
		\begin{equation*}
			\{H(i\zeta_{j}),H(i\zeta_{-j})\}_{j=1}^{\nu_{p}}\quad\text{and}\quad\{H(i\omega_{k}),H(i\omega_{-k})\}_{k=1}^{\nu_{q}}.
		\end{equation*}
The matrices $\widetilde{\mathbb M}, \widetilde{\mathbb D}, \widetilde{\mathbb K}$ 
can be partitioned into $2\times 2$ blocks, $\widetilde{\mathbb M}^{(2)}_{k,j}, 
\widetilde{\mathbb D}^{(2)}_{k,j}, \widetilde{\mathbb K}^{(2)}_{k,j}$, which are 
compatible with the conjugate pairs for the quadrature rule, along with the nodes 
$(i\omega_k, -i\omega_k)$ and $(i\zeta_j, -i\zeta_k)$. It follows  from 
(\ref{relation_damp}) that 
\begin{equation*}
	\left[
	\begin{array}{cc}
		i\omega_k&\\
		&-i\omega_k
	\end{array}
	\right]^2\widetilde{\mathbb{M}}_{k,j}^{(2)}+\left[
	\begin{array}{cc}
		i\omega_k&\\
		&-i\omega_k
	\end{array}
	\right]\left(\alpha\widetilde{\mathbb{M}}_{k,j}^{(2)}+\beta\widetilde{\mathbb{K}}_{k,j}^{(2)}\right)+
	\widetilde{\mathbb{K}}_{k,j}^{(2)}=\mu_k\rho_j\zeta_j\left[
	\begin{array}{c}
		1\\
		1
	\end{array}
	\right]\left[H(i\zeta_j) \quad -\overline{H(i\zeta_j)}\right],
\end{equation*}
\begin{equation*}
	\widetilde{\mathbb{M}}_{k,j}^{(2)}\left[
	\begin{array}{cc}
		i\zeta_j&\\
		&-i\zeta_j
	\end{array}
	\right]^2+\left(\alpha\widetilde{\mathbb{M}}_{k,j}^{(2)}+\beta\widetilde{\mathbb{K}}_{k,j}^{(2)}\right)
	\left[
	\begin{array}{cc}
		i\zeta_j&\\
		&-i\zeta_j
	\end{array}
	\right]+
	\widetilde{\mathbb{K}}_{k,j}^{(2)}=\mu_k\rho_j\zeta_j\left[
	\begin{array}{c}
		H(i\omega_k)\\
		\overline{H(i\omega_k)}
	\end{array}
	\right]\left[1 \quad -1\right].
\end{equation*}
We define the unitary matrices 
		  \begin{equation*}
		  	\mathcal{F}_1 = \dfrac{1}{\sqrt{2}}\left[\begin{array}{cc}
		  		1	& -1 \\
		  		i	&  i
		  	\end{array}\right],\quad \mathcal{F}_2 = 
		  	\dfrac{1}{\sqrt{2}}\left[\begin{array}{cc}
		  	1	& 1 \\
		  	i	& -i
	  	\end{array}\right],
		  \end{equation*}
and the following transformation 
\begin{equation}\label{transformation}
	\widetilde{\mathbb{M}}_{k,j}^{(2)}=\mathcal{F}_2^{\mathrm 
	H}\widetilde{\mathbb{M}}_{k,j}^{\mathrm R}\mathcal{F}_1, 
\widetilde{\mathbb{K}}_{k,j}^{(2)}=\mathcal{F}_2^{\mathrm 
	H}\widetilde{\mathbb{K}}_{k,j}^{\mathrm R}\mathcal{F}_1. 
\end{equation}
Substituting (\ref{transformation}) into the above equations and performing some 
matrix operation lead to 
\begin{equation}\label{real_equa}
	\begin{split}
		-\Omega_{k,2}\widetilde{\mathbb{M}}_{k,j}^{\mathrm 
		R}+\Omega_{k,1}\left(\alpha\widetilde{\mathbb{M}}_{k,j}^{\mathrm 
	R}+\beta\widetilde{\mathbb{K}}_{k,j}^{\mathrm R}\right)+
	\widetilde{\mathbb{K}}_{k,j}^{\mathrm R}&=H_{\zeta, kj},\\
	-\widetilde{\mathbb{M}}_{k,j}^{\mathrm 
	R}\Theta_{j,2}+\left(\alpha\widetilde{\mathbb{M}}_{k,j}^{\mathrm 
		R}+\beta\widetilde{\mathbb{K}}_{k,j}^{\mathrm R}\right)\Theta_{j,1}+
	\widetilde{\mathbb{K}}_{k,j}^{\mathrm R}&=H_{\omega,kj},
	\end{split}
\end{equation}
where the coefficient matrices are defined as 
\begin{equation*}
	\begin{split}
	&\Omega_{k,2}=\left[
	\begin{array}{cc}
		\omega_k^2&\\
		&\omega_k^2
	\end{array}
	\right], \Omega_{k,1}=\left[
	\begin{array}{cc}
		&\omega_k\\
		-\omega_k&
	\end{array}
	\right], H_{\zeta, kj}=2\mu_k\rho_j\zeta_j\left[
	\begin{array}{cc}
		\mathrm{Re}(H(i\zeta_j))&\mathrm{Im}(H(i\zeta_j))\\
		0&0
	\end{array}
	\right],\\
	&\Theta_{j,2}=\left[
	\begin{array}{cc}
		\zeta_j^2&\\
		&\zeta_j^2
	\end{array}
	\right], \Theta_{j,1}=\left[
	\begin{array}{cc}
		&\zeta_j\\
		-\zeta_j&
	\end{array}
	\right], H_{\omega,kj}=2\mu_k\rho_j\zeta_j\left[
	\begin{array}{cc}
		\mathrm{Re}(H(i\omega_k))&0\\
		-\mathrm{Im}(H(i\omega_k))&0
	\end{array}
	\right]. 
	\end{split}
\end{equation*}
Consequently, we can obtain $\widetilde{\mathbb{M}}_{k,j}^{\mathrm R}$ and 
$\widetilde{\mathbb{K}}_{k,j}^{\mathrm R}$ by solving linear equations 
(\ref{real_equa}) in real arithmetic. Similarly, we partition 
$\mathbb{\widetilde{\mathbb B}}$ and 
$\mathbb{\widetilde{\mathbb C}}$ into $2\times 1$ and $1\times 2$ blocks, 
respectively, and define the 
transformation 
$$\widetilde{\mathbb B}_k^{(2)}=\mathcal{F}_2^{\mathrm 
H}\widetilde{\mathbb B}_k^{\mathrm R} \quad \text{and} \quad \widetilde
{\mathbb C}_j^{(2)}=\widetilde{\mathbb C}_j^{\mathrm R}\mathcal{F}_1.$$
One can check directly that $\widetilde{\mathbb B}_k^{\mathrm R}, 
\widetilde{\mathbb C}_j^{\mathrm R}$ are real-value vectors as follows
\begin{equation}\label{eqn_B_C_R}
	\widetilde{\mathbb B}_k^{\mathrm R}=\sqrt{2}\mu_k
	\left[
	\begin{array}{c}
		\mathrm{Re}(H(i\omega_k))\\
		-\mathrm{Im}(H(i\omega_k))
	\end{array}
	\right], 
	\widetilde{\mathbb C}_j^{\mathrm R}=\sqrt{2}\rho_j\zeta_j
	\left[\mathrm{Re}(H(i\zeta_j))\quad \mathrm{Im}(H(i\zeta_j)) 
	\right].
\end{equation}
Then the real-valued counterparts corresponding to 
$\widetilde{\mathbb{M}}$, 
$\widetilde{\mathbb{D}}$, $\widetilde{\mathbb{K}}$, $\widetilde{\mathbb{B}}$, and 
$\widetilde{\mathbb{C}}$ are given by
		  \begin{eqnarray*}
		  	&&\widetilde{\mathbb{M}}^R=(I_{\nu_q}\otimes 
		  	\mathcal{F}_2)\widetilde{\mathbb{M}}(I_{\nu_p}\otimes 
		  	\mathcal{F}_1^\mathrm{H}),\\
		  	&&\widetilde{\mathbb{D}}^R=(I_{\nu_q}\otimes 
		  	\mathcal{F}_2)\widetilde{\mathbb{D}}(I_{\nu_p}\otimes 
		  	\mathcal{F}_1^\mathrm{H}),\\
		  	&&\widetilde{\mathbb{K}}^R=(I_{\nu_q}\otimes 
		  	\mathcal{F}_2)\widetilde{\mathbb{K}}(I_{\nu_p}\otimes 
		  	\mathcal{F}_1^\mathrm{H}),\\
		  	&&\widetilde{\mathbb{B}}^R=(I_{\nu_q}\otimes 
		  	\mathcal{F}_2)\widetilde{\mathbb{B}},\\
		  	&&\widetilde{\mathbb{C}}^R=\widetilde{\mathbb{C}}(I_{\nu_p}\otimes 
		  	\mathcal{F}_1^\mathrm{H}),
		  \end{eqnarray*}
where $\otimes$ denotes the Kronecker product. In practice, one can replace the 
$\widetilde{\mathbb M}$, $\widetilde{\mathbb D}$, $\widetilde{\mathbb K}$, 
$\widetilde{\mathbb B}$ and $\widetilde{\mathbb C}$  with the counterparts 
$\widetilde{\mathbb{M}}^R, \widetilde{\mathbb{D}}^R, \widetilde{\mathbb{K}}^R, 
\widetilde{\mathbb{B}}^R, \widetilde{\mathbb{C}}^R$ in Algorithm 
\ref{al:dataSOBT}, which avoids the complex arithmetic and results in real-value 
reduced models. 

\begin{remark}
	In \cite{Schulze2018,Gosea2024}, the structured realizations of second-order 
	systems is considered which impose the interpolation conditions at the 
	interpolation points $\lambda_i$ and $\mu_i$ for $i\in\{1, \cdots, l\}$. With 
	the notations
	\begin{equation*}
		\begin{split}
			\mathbb{I}\in\mathbb{R}^{l\times 1}&=[1, 1, \cdots, 1]^{\top},\\
			\Lambda=\mathrm{diag}\{\lambda_1, \cdots, \lambda_l\}, &
			\Psi=\mathrm{diag}\{\mu_1, \cdots, \mu_l\},\\
			\hat{H}(\Lambda)=[H(\lambda_1), \cdots, H(\lambda_l)]^{\top}, &
			\hat{H}(\Psi)=[H(\mu_1), \cdots, H(\mu_l)]^{\top},
		\end{split}
	\end{equation*}  
the coefficient matrices $M_\ell, K_\ell, D_\ell$ of reduced models are required to 
satisfy the equalities
\begin{equation*}
\begin{split}
	&M_\ell\Lambda^2+D_\ell\Lambda+K_\ell=\hat{H}(\Psi)\mathbb{I}^{\top},\\
	&\Psi^2M_\ell+\Psi D_\ell+K_\ell=\mathbb{I}\hat{H}(\Lambda)^{\top},
\end{split}
\end{equation*}	
which are almost the same as (\ref{eqMDK1}) and (\ref{eqMDK2}) but with different 
weights of the frequency response. However, the resulting reduced models produced 
by Algorithm \ref{al:dataSOBT} have no any interpolation property in general. 		
\end{remark}

\subsection{Low-rank execution based on Sylvester equations}
The focus of the data-driven BT provided in Algorithm \ref{al:dataSOBT} is to 
approximate the main qualities via the frequency response of systems, where the 
number of quadrature nodes $N_p, N_q$ effects greatly the accuracy of the 
approximation. A large number of quadrature nodes results in high dimensional 
matrices  $\widetilde{\mathbb{M}}, \widetilde{\mathbb{D}}, \widetilde{\mathbb{K}}$. 
Although the computation load can be reduced dramatically based on 
$\widetilde{\mathbb{M}}^{\mathrm R}, \widetilde{\mathbb{D}}^{\mathrm R}, 
\widetilde{\mathbb{K}}^{\mathrm R}$ in real arithmetic,  
and the element-wise assembly of these matrices is still unacceptably time and 
storage consuming in practice. Besides, the SVD of $\widetilde{\mathbb{M}}^{\mathrm 
R}$ in step 
2 of Algorithm $\ref{al:dataSOBT}$ costs too much in the large scale settings, 
making the whole data-driven BT procedure relatively inefficiency. Inspired by the 
techniques in \cite{Hamadi2023}, we 
switch to Sylvester equations the main 
quantities in Algorithm \ref{al:dataSOBT} satisfy, and calculate the low-rank 
approximation to these matrices in order to reduce the computational load 
dramatically in practice.

Consider linear equations (\ref{real_equa}). The real-value matrices 
$\widetilde{\mathbb{M}}^{\mathrm R}$ and $\widetilde{\mathbb{K}}^{\mathrm R}$ 
satisfy 
\begin{equation*}
	\begin{split}
	-\Omega_2 \widetilde{\mathbb{M}}^{\mathrm R}+\Omega_1(\alpha 
	\widetilde{\mathbb{M}}^{\mathrm R}+\beta\widetilde{\mathbb{K}}^{\mathrm 
	R})+\widetilde{\mathbb{K}}^{\mathrm 
	R}=H_{\zeta},\\
	-\widetilde{\mathbb{M}}^{\mathrm R}\Theta_2 +(\alpha 
	\widetilde{\mathbb{M}}^{\mathrm R}+\beta\widetilde{\mathbb{K}}^{\mathrm 
	R})\Theta_1+\widetilde{\mathbb{K}}^{\mathrm 
	R}=H_{\omega},
	\end{split}
\end{equation*}
where the coefficient matrices are defined as
\begin{equation*}
	\begin{split}
	H_{\zeta}=\{H_{\zeta, kj}\}_{kj}\in\mathbb{R}^{2\nu_q\times 2\nu_p}, 
	H_{\omega}=\{H_{\omega, kj}\}_{kj}\in\mathbb{R}^{2\nu_q\times 2\nu_p},\\	
	\Omega_2=\mathrm{diag}(\Omega_{1,2}, \cdots, \Omega_{\nu_q,2}),
	\Omega_1=\mathrm{diag}(\Omega_{1,1}, \cdots, \Omega_{\nu_q,1}),\\
	\Theta_2=\mathrm{diag}(\Theta_{1,2}, \cdots, \Theta_{\nu_p,2}),
	\Theta_1=\mathrm{diag}(\Theta_{1,1}, \cdots, \Theta_{\nu_p,1}).\\
	\end{split}
\end{equation*}
The direct algebraic manipulations lead to 
\begin{equation}
	\begin{split}
		&(\beta\Omega_1+I)^{-1}(-\Omega_2+\alpha\Omega_1)\widetilde{\mathbb 
		M}^{\mathrm R}-\widetilde{\mathbb M}^{\mathrm 
		R}(-\Theta_2+\alpha\Theta_1)(\beta\Theta_1+I)^{-1}\\ 
	    =&(\beta\Omega_1+I)^{-1}H_{\zeta}-
		H_{\omega}(\beta\Theta_1+I)^{-1},\label{sylv_eqn_1}
	\end{split}
\end{equation}
\begin{equation}
	\begin{split}
		&(-\Omega_2+\alpha\Omega_1)^{-1}(\beta\Omega_1+I)\widetilde{\mathbb 
		K}^{\mathrm
		 R}-\widetilde{\mathbb K}^{\mathrm 
		R}(\beta\Theta_1+I)(-\Theta_2+\alpha\Theta_1)^{-1}\\ 
	=&(-\Omega_2+\alpha\Omega_1)^{-1}H_{\zeta}-H_{\omega}(-\Theta_2+\alpha\Theta_1)^{-1},\label{sylv_eqn_2}
	\end{split}
\end{equation}
which implies that $\widetilde{\mathbb{M}}^{\mathrm R}$ and 
$\widetilde{\mathbb{K}}^{\mathrm R}$ can be obtained by solving Sylvester 
equations independently.  
Note that the coefficient matrices 
$$\beta\Omega_1+I, \beta\Theta_1+I, -\Omega_2+\alpha\Omega_1, 
-\Theta_2+\alpha\Theta_1$$ are block diagonal 
matrices along with $2\times 2$ diagonal blocks, and one can build up  
(\ref{sylv_eqn_1}) and (\ref{sylv_eqn_2}) at far lower cost. Furthermore, when 
$\alpha, \beta\ge 
0$, $\omega_k, \zeta_j>0$ for $k=1, 2, \cdots, \nu_q, j=1, 2, \cdots, \nu_p$, they 
are all invertible matrices. Further, because of 
$\mathrm{rank}(H_{\omega})=1$ and $\mathrm{rank}(H_{\zeta})=1$, generally there hold 
\begin{equation*}
	\begin{split}
	\mathrm{rank}((\beta\Omega_1+I)^{-1}H_{\zeta}-
	H_{\omega}(\beta\Theta_1+I)^{-1})=2, \\
	\mathrm{rank}((-\Omega_2+\alpha\Omega_1)^{-1}H_{\zeta}-H_{\omega}
	(-\Theta_2+\alpha\Theta_1)^{-1})=2. 
	\end{split}
\end{equation*}	
In fact, with the following notations 
$$\Delta_l=\mathrm{diag}(2\mu_1, 2\mu_2, \cdots, 
2\mu_{\nu_q})\otimes I_2, \Delta_r=\mathrm{diag}(\rho_1\zeta_1, \rho_2\zeta_2, 
\cdots, \rho_{\nu_p}\zeta_{\nu_p})\otimes I_2,$$
$$l_{\zeta}=[\mathrm{Re}(H(i\zeta_1)) \quad\mathrm{Im}(H(i\zeta_1))  \quad\cdots 
\quad 
\mathrm{Re}(H(i\zeta_{\nu_p})) \quad\mathrm{Im}(H(i\zeta_{\nu_p}))]^{\top},$$
$$l_{\omega}=[\mathrm{Re}(H(i\omega_1)) \quad-\mathrm{Im}(H(i\omega_1))  
\quad\cdots \quad 
\mathrm{Re}(H(i\omega_{\nu_q})) \quad-\mathrm{Im}(H(i\omega_{\nu_q}))]^{\top},$$
$$e_{\omega}=[1 \quad0 \quad \cdots \quad 1 \quad
0]^{\top}\in\mathbb R^{2\nu_q\times 1},  e_{\zeta}=[1 \quad
0 \quad\cdots\quad 1 \quad0]^{\top}\in\mathbb R^{2\nu_p \times 1},$$
the right sides of (\ref{sylv_eqn_1}) and (\ref{sylv_eqn_2}) have the factor form 
\begin{equation*}
	\begin{split}
(\beta\Omega_1+I)^{-1}H_{\zeta}-H_{\omega}(\beta\Theta_1+I)^{-1}&=
\begin{bmatrix}
(\beta\Omega_1+I)^{-1}\Delta_l e_{\omega} & \Delta_ll_{\omega}
\end{bmatrix}
	\begin{bmatrix}
	l_{\zeta}^{\top}\Delta_r\\
	-e_{\zeta}^{\top}\Delta_r(\beta\Theta_1+I)^{-1}
	\end{bmatrix},\\
(-\Omega_2+\alpha\Omega_1)^{-1}H_{\zeta}-H_{\omega}(-\Theta_2+\alpha\Theta_1)^{-1}&=
\begin{bmatrix}
	(-\Omega_2+\alpha\Omega_1)^{-1}\Delta_l e_{\omega} & \Delta_ll_{\omega}
\end{bmatrix}
\begin{bmatrix}
	l_{\zeta}^{\top}\Delta_r\\
	-e_{\zeta}^{\top}\Delta_r(-\Theta_2+\alpha\Theta_1)^{-1}
\end{bmatrix},
\end{split}
\end{equation*}
respectively. As a consequent, one can expect the accurate low-rank approximate 
solution of (\ref{sylv_eqn_1}) and (\ref{sylv_eqn_2}) based on the Krylov subspace 
methods.

For ease of presentation, we adopt the following compact formulation of 
(\ref{sylv_eqn_1}) to show the main procedure of the computation
\begin{equation}\label{eq:sylvester}
	Z\widetilde{\mathbb M}^{\mathrm R}-\widetilde{\mathbb M}^{\mathrm R}Y=EF^T,
\end{equation}
where $Z \in \mathbb{R}^{2\nu_q\times 2\nu_q}$, $Y \in \mathbb{R}^{2\nu_p \times 
2\nu_p}$, $E 
\in \mathbb{R}^{2\nu_q \times 2}$, $F \in \mathbb{R}^{2\nu_p \times 2}$ are the 
counterparts of coefficient matrices of (\ref{sylv_eqn_1}). Below we 
employ the extended Krylov subspace method to solve (\ref{eq:sylvester}) 
approximately.
The extended block Krylov subspace is defined as
\begin{equation}
	\mathcal{K}^{\mathrm 
	{ext}}_m(G,J)=\mathcal{K}_m(G,J)\cup\mathcal{K}_m(G^{-1},G^{-1}J),
\end{equation}
where $\mathcal{K}_m(G,J)=\mathrm{span}(J, GJ, \cdots, G^{m-1}J)$ is the standard 
Krylov subspace for the given matrices $G\in \mathbb{R}^{N\times N},J \in 
\mathbb{R}^{N \times s}$  \cite{Simoncini2007}. An orthogonal basis of the extended 
Krylov subspace can be 
extracted by using the the extended
block Arnoldi process. The basic idea of the extended Krylov subspace method is to
restrict the solution of (\ref{eq:sylvester}) onto the extended Krylov subspaces 
spanned by its coefficient matrices. We apply simultaneously the extended
block Arnoldi process to the pairs $\mathcal{K}^{\mathrm 
{ext}}_m(Z,E)$ and $\mathcal{K}^{\mathrm 
{ext}}_m(Y^{\top},F)$. After $m$ iterations, we get two sets of orthogonal basis 
$\mathbb{V}^Z_m 
 \in \mathbb{R}^{2\nu_q \times 4m}$ and $\mathbb{V}^Y_m \in \mathbb{R}^{2\nu_p 
 \times 4m}$ associated with the two extended Krylov subspaces in general, 
 respectively. The approximate 
 solution is given as
\begin{equation}\label{app_solut}
	X_m = \mathbb{V}^Z_mS_m(\mathbb{V}^Y_m)^\mathrm{\top},
\end{equation}
where $S_m \in \mathbb{R}^{4m \times 4m}$ is a small matrix to be determined. 
Substituting $X_m$ into (\ref{eq:sylvester}) leads to the residual 
$R_m = ZX_m-X_mY-EF^\mathrm{\top}. $
By enforcing the Galerkin condition there holds
\begin{equation*}
	(\mathbb{V}^Z_m)^\mathrm{\top}R_m\mathbb{V}^Y_m=0. 
\end{equation*}
As both $\mathbb{V}^Z_m$ and $\mathbb{V}^Y_m$ are the orthonormal basis, the 
Galerkin condition 
boils down to the following linear equation about $S_m$ 
\begin{equation*}
	\begin{split}
	0&=(\mathbb{V}^Z_m)^{\top}(ZX_m-X_mY-EF^{\top})\mathbb{V}^Y_m\\
	&=(\mathbb{V}^Z_m)^{\top}Z\mathbb{V}^Z_mS_m-S_m(\mathbb{V}^Y_m)^{\top}Y
	\mathbb{V}^Y_m-(\mathbb{V}^Z_m)^{\top}EF^{\top}\mathbb{V}^Y_m. 
	\end{split}
\end{equation*}
With the notations $\mathbb{T}^Z_m=(\mathbb{V}^Z_m)^{\top}Z\mathbb{V}^Z_m$, 
$\mathbb{T}^Y_m=(\mathbb{V}^Y_m)^{\top}Y
\mathbb{V}^Y_m$, $E_m = (\mathbb{V}^Z_m)^{\top}E$, $F_m = 
(\mathbb{V}^Y_m)^{\top}F$, we get the following low-dimensional Sylvester 
equation
\begin{equation}\label{low:sylvesyer}
	\mathbb{T}^Z_mS_m-S_m\mathbb{T}^Y_m=E_mF^{\top}_m. 
\end{equation}
As a result, a low-rank approximation to $\widetilde{\mathbb M}^{\mathrm R}\approx 
X_m$ is obtained by solving (\ref{low:sylvesyer}) via a direct method as described 
in \cite{Simoncini2016}.

Likewise, we can also get the low-rank approximation to $\widetilde{\mathbb 
K}^{\mathrm R}$. Note that our main purpose is to enable an efficient execution of 
Algorithm \ref{al:dataSOBT}. From this point of view, there is no need to assemble 
$X_m$ explicitly, and the factor form of (\ref{app_solut}) is more preferable. We 
need the SVD of $\widetilde{\mathbb M}^{\mathrm R}\approx X_m$ in step 2 of 
Algorithm \ref{al:dataSOBT}. By computing the SVD of low order matrix $S_m$, which 
is referred as $[U, \Sigma, W]=\mathrm{SVD}(S_m)$, the SVD of $X_m$ reads
$$\widetilde{\mathbb M}^{\mathrm R}\approx X_m=\mathbb{V}^Z_mU\Sigma 
(\mathbb{V}^Y_mW)^\mathrm{\top},$$
which benefits the execution of the proposed approach a lot. We summarize 
the main steps of the low-rank execution of the proposed approach in real arithmetic 
in Algorithm \ref{alg:dbt-k}.

\begin{algorithm}[htb]
	\caption{Data-driven BT of second order systems in real arithmetic along with 
	low-rank approximation. }
	\label{alg:dbt-k}
	\begin{algorithmic}[1]
		\Require
		Quadrature nodes $\zeta_j, \omega_k>0$ and weights $\rho_j, 
		\varphi_k$, for $j = 1\dots \nu_p, k = 1\dots 
		\nu_q$;
		
		Sample data of $H(s)$ at the quadrature nodes, and the 
		index $1 \leq r \leq \mathrm{min}(N_p, N_q)$.
		
		\Ensure
		Reduced models $M_r\in\mathbb{R}^{r\times 
			r},D_r\in\mathbb{R}^{r\times r},K_r\in\mathbb{R}^{r\times 
			r},B_r\in\mathbb{R}^{r\times 1},C_r\in\mathbb{R}^{1\times r}$. 
		\State Use the sample data to assemble the coefficient matrices $Z, Y, E, F$ 
		of Sylvester equation (\ref{eq:sylvester}).
		
		\State Apply the extended block Arnoild process to the pair 
		$\mathcal{K}^{\mathrm 
			{ext}}_m(Z,E)$ and $\mathcal{K}^{\mathrm 
			{ext}}_m(Y^{\top},F)$, and produce the orthogonal basis 
		$\mathbb{V}^Z_m $ and $\mathbb{V}^Y_m$. 
		
		\State Construct the low-order Sylvester equation (\ref{low:sylvesyer}) and 
		solve it via the standard methods directly. 
		
		\State Compute the SVD of the matrix $S_m$
		\begin{equation}
			S_m=\left[\begin{array}{cc}
				U_1	& U_2
			\end{array}\right]\left[\begin{array}{cc}
				\Sigma_1&  \\
				& \Sigma_2
			\end{array}\right]\left[\begin{array}{c}
				W^{\top}_1\\
				W^{\top}_2
			\end{array}\right],
		\end{equation}
		where $\Sigma_1\in \mathbb{R}^{r\times r}$ and $\Sigma_2\in \mathbb{R}^{(4m-r)\times(4m-r)}$.
		
		\State Compute the low-rank approximation $\widetilde{\mathbb K}^{\mathrm 
		R}$ by solving the low-order Sylvester equation corresponding to 
		(\ref{sylv_eqn_2}). Compute $\widetilde{\mathbb 
			B}^{\mathrm R}$ and $\widetilde{\mathbb C}^{\mathrm R}$ via 
		(\ref{eqn_B_C_R}). 
		
		\State The reduced models are given by
		\begin{equation*}
			\begin{array}{c}
				M_r=I_r,\quad
				K_r=\Sigma_{1}^{-1 / 
				2}U^{\top}_1(\mathbb{V}_m^{Z})^{\top}\widetilde{\mathbb K}^{\mathrm 
					R}\mathbb{V}_m^{Y}W_1\Sigma_{1}^{-1 / 2},\quad	
				D_r=\alpha M_r+\beta K_r,	\\	
				B_r=\Sigma_{1}^{-1 / 2}
			U^{\top}_1(\mathbb{V}_m^{Z})^{\top}\widetilde{\mathbb B}^{\mathrm 
					R},\quad
				C_r=\widetilde{\mathbb C}^{\mathrm 
					R}\mathbb{V}_m^{Y}W_1\Sigma_{1}^{-1 / 2}.
			\end{array}
		\end{equation*}
	\end{algorithmic}
\end{algorithm}

\begin{remark}
	The coefficient matrices $\mathbb{T}^Z_m$ and $\mathbb{T}^Y_m$ are the 
	restrictions of $Z$ and $Y$ to the extended Krylov subspace. Both are block 
	upper Hessenberg matrices, and can be derived efficiently via the extended block 
	Arnoild procedure, instead of any extra inner products of long vectors. We use 
	the same number of 
	iterates for the pair extended Krylov subspace in the above. In practice, one 
	can 
	determine the number of iterates individually for each subspace in light of the 
	residual in a certain norm computed in an efficient way \cite{Hamadi2023}. 
\end{remark}

\begin{remark}
	In the above we focus on the SISO case, and the MIMO case can be immediately 
	developed similarly. The sample of $H(s)$ at a quadrature node is a $p\times m$ 
	matrix in MIMO case. The $(k, j)$ elements of 
	$\widetilde{\mathbb M}$, $\widetilde{\mathbb K}$ and $\widetilde{\mathbb D}$ are 
	matrices of order $p\times m$ in the sense of block matrix, and the linear 
	equations similar to the ones in Proposition 1 can be defined using the same 
	strategy. So, Algorithm \ref{al:dataSOBT} proceeds naturally with the new 
	quantities. Moreover, with the symmetic nodes along the real axis, the unitary 
	matrices 
	\begin{equation*}
		\mathcal{F}_1 = \dfrac{1}{\sqrt{2}}\left[\begin{array}{cc}
			1	& -1 \\
			i	&  i
		\end{array}\right]\otimes I_m,\quad \mathcal{F}_2 = 
		\dfrac{1}{\sqrt{2}}\left[\begin{array}{cc}
			1	& 1 \\
			i	& -i
		\end{array}\right]\otimes I_p
	\end{equation*}
	can be use to  define a proper transformation, which enables an execution in 
	real arithmetic. The low-rank approximation is also valid in MIMO case. We omit 
	the details in this paper for brevity. 
\end{remark}

\begin{remark}
	Generally, there are four types of balanced form for second order systems, 
	namely, velocity balanced, position balanced, position-velocity balanced, and 
	velocity-position balanced.  We take the velocity balanced truncation as an 
	example for our discussion. In fact, the data-driven BT presented in this paper 
	can also 
	be executed in the framework of position-velocity balanced form, with some minor 
	modification. However, because of the 
	relative complex expression of $Q_p$ in (\ref{G-3}), the proposed approach is 
	not available in the case of position balanced and velocity-position balanced 
	form. As for the data-driven BT for the general second order systems (no 
	explicit relationship among coefficient matrices), these results warrant further 
	investigation. 
\end{remark}

\section{Numerial examples}\label{sec:sec-4}

In this section, we use three numerical examples to illustrate the effectiveness of 
our approach. The simulation is performed via Matlab (R2023a) on a laptop with 
Intel(R) Core(TM) Ultra 5 125H and 16 GB RAM.

We execute Algorithm \ref{al:dataSOBT} 
and Algorithm \ref{alg:dbt-k} to produce the reduced models Data-BT-SOPD and 
KryData-BT-SOPD, respectively. Although Algorithm \ref{al:dataSOBT} is presented in 
complex arithmetic in subsection 3.2, we carry out its counterpart version in real 
arithmetic in the simulation to generate reduced models with real-value 
coefficient matrices. In Algorithm \ref{alg:dbt-k}, the related Sylvester equations 
are projected via the extended Arnoldi procedure first, and the projected low 
dimensional Sylvester equations are computed via the Matlab function lyap. We simply 
use the same number of nodes for the quadrature approximation to $P_v$ and $Q_v$. 
The Matlab function logspace is employed to get the logarithmically spaced points 
$\omega_k$ and $\zeta_j$, and the transfer functions are evaluated at 
$\omega_k\mathrm {i}, \zeta_j\mathrm {i}$ to provide the measurements. 
The 
quadrature weights are determined according to the trapezoid 
quadrature rule in the simulation. 
We refer the reader to \cite{gosea2022} for more details on other representative 
quadrature rules. The 
standard BT of second order systems shown in 
Algorithm \ref{SOBT} is also implemented for the purpose of comparison, and the 
associated reduced models are referred as BT-SOPD.

\begin{example}
	We consider the building model from the SLICOT library \cite{chahlaoui2002}. It 
	describes a 
	building with 8 floors each having 3 degrees of freedom, and has 24 variables in 
	the form of the second order systems. The Rayleigh damping coefficients in the 
	simulation are $\alpha=0.05$ and $\beta=0.05$. 
	
	\begin{figure}[htb]
		\centering
		\begin{minipage}{0.4\textwidth} 
			\centering 
			\includegraphics[width=0.9\textwidth]{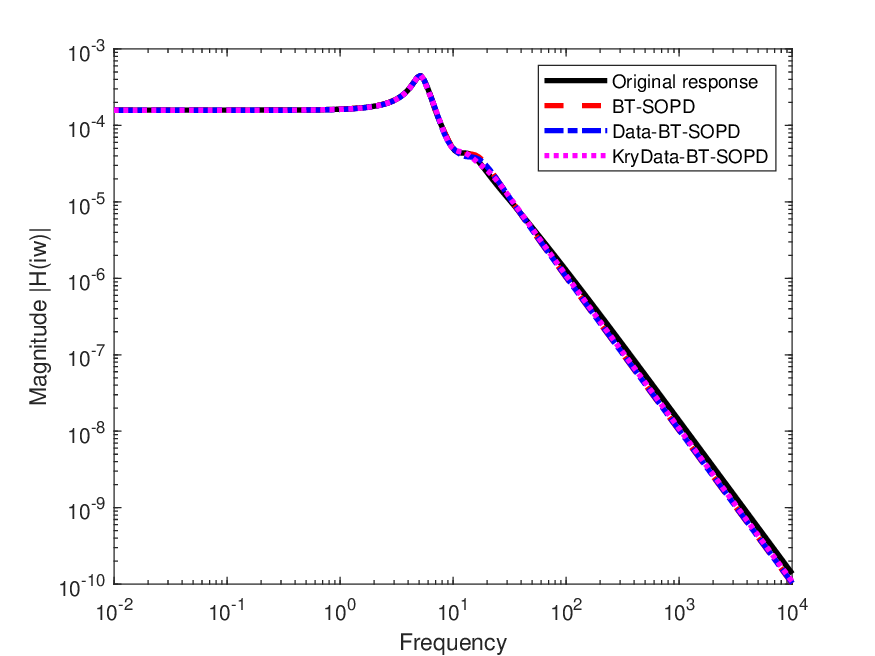}
		\end{minipage}
		\begin{minipage}{0.4\textwidth} 
			\centering 
			\includegraphics[width=0.9\textwidth]{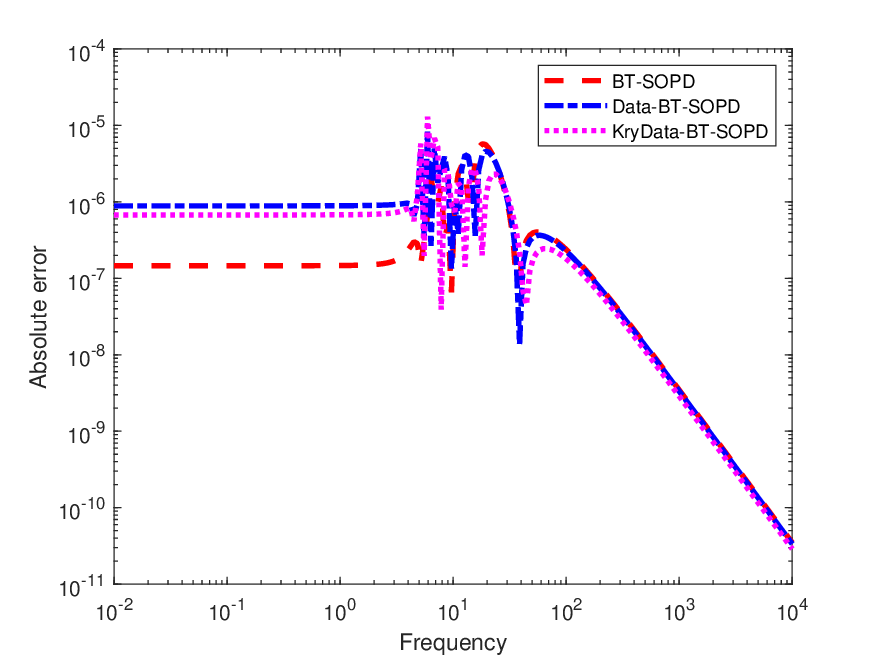}
		\end{minipage}
		\caption{The frequency response (left) and the absolute errors (right) with 
			$r=4$, $N_q, N_p=50$ and $m=5$.}
		\label{fig11}
	\end{figure}

\begin{figure}[htb]
	\centering
	\begin{minipage}{0.4\textwidth} 
		\centering 
		\includegraphics[width=0.9\textwidth]{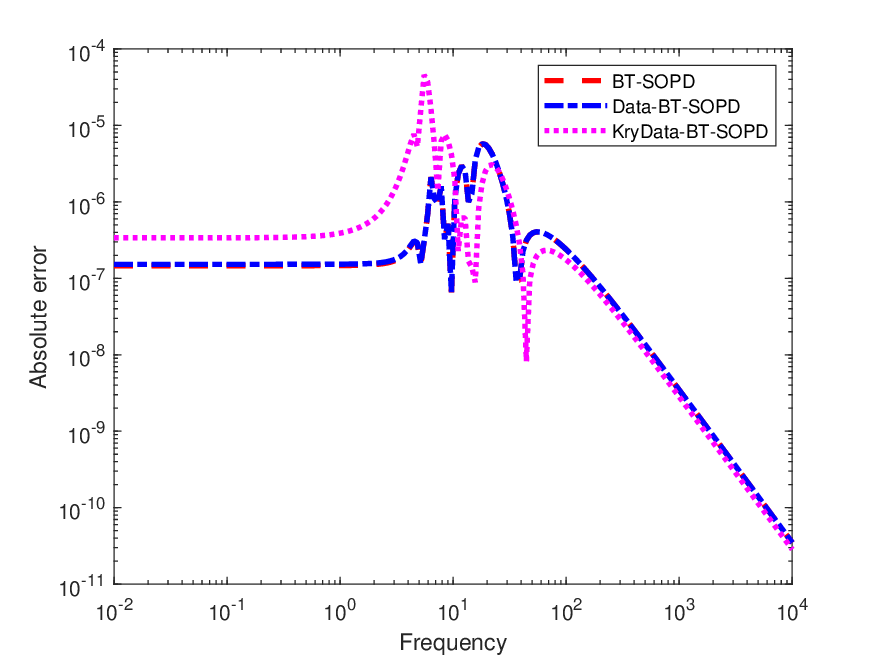}
	\end{minipage}
	\begin{minipage}{0.4\textwidth} 
		\centering 
		\includegraphics[width=0.9\textwidth]{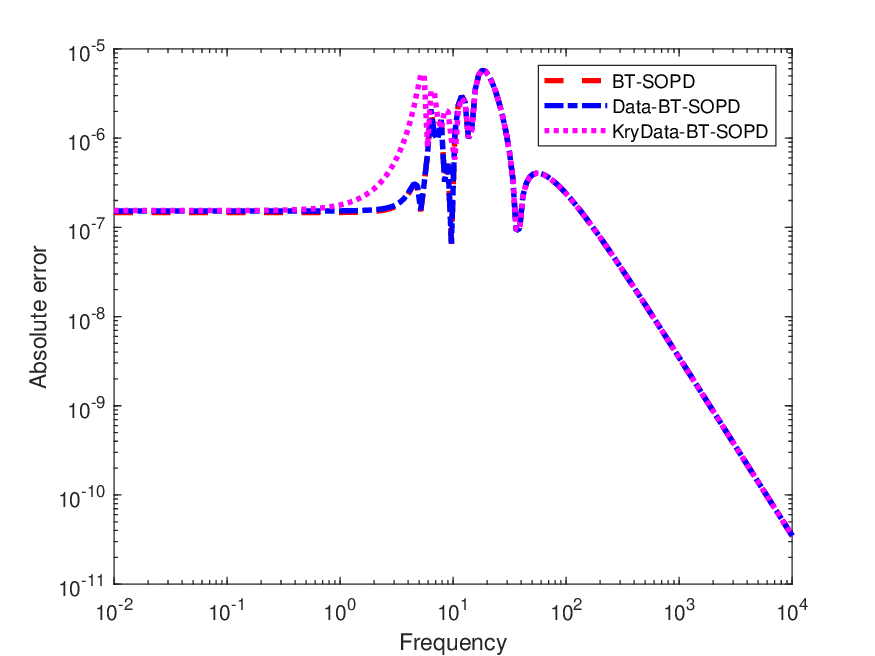}
	\end{minipage}
	\caption{The absolute errors with $r=4, N_q, N_p=500, m=5$ (left) and $r=4, 
		N_q, N_p=500, m=20$ (right).}
	\label{fig12}
\end{figure}
	
	We first select 50 nodes in the interval $s_i\in[10^{-2}, 
	10^{4}]$ logarithmically for the numerical quadrature. The 
	nodes and the 
	conjugations are 
	partitioned into two parts, and then reordered as shown in (\ref{order-nodes}) 
	in the simulation. The evaluation of $H(\mathrm{i} s_i)$ at these nodes are 
	collected to produce 100 measurements. With the reduced order $r=4$, three 
	reduced models are generated by implementing Algorithm 
	\ref{SOBT}-\ref{alg:dbt-k}. Note that we set $m=5$ in the execution of Algorithm 
	\ref{alg:dbt-k}. Figure \ref{fig11} depicts the frequency domain 
	responses and the associated absolute errors for each reduced model. The 
	original system is well approximated by all reduced models, and we can not 
	distinguish them clearly from the response depiction. However, one can observe 
	an evident distinction between the standard BT and the proposed data-driven 
	version from the error depiction, especially in the lower 
	frequency domain. This is due to the relatively less samples involved in the 
	quadrature.

	We then enrich the quadrature nodes and use the relatively high order extended 
	Krylov subspace in our methods. Figure \ref{fig12} shows the evolution of 
	absolute errors with respect to different choices of $N_q, N_p$ and $m$. The 
	behavior of reduced models produced via our methods approaches gradually to that 
	of the standard BT in the simulation. By taking 500 nodes in the same interval, 
	DataBT-SOPD  mimics accurately BT-SOPD, and the deviation between Data-BT-SOPD 
	and KryData-BT-SOPD declines dramatically for $m=20$. In fact, when $m=30$ is 
	adopted for this example, one can hardly differentiate each reduced model from 
	the error depiction. Given the input function $u(t)=\mathrm{e}^{-t}\sin(t)$, 
	Figure 
	\ref{fig13} plots the outputs of reduced models in the time domain. The main 
	purpose of Algorithm \ref{alg:dbt-k} is to provide a good 
	approximation and thereby to enable an efficient execution of the data-driven 
	BT approach. The CPU time spent in the construction of reduced models is listed 
	in Table \ref{tab:11}, where the advantage of Algorithm \ref{alg:dbt-k} is 
	exhibited evidently in the intensive sample data.

\begin{figure}[htb]
	\centering
	\begin{minipage}{0.4\textwidth} 
		\centering 
		\includegraphics[width=0.9\textwidth]{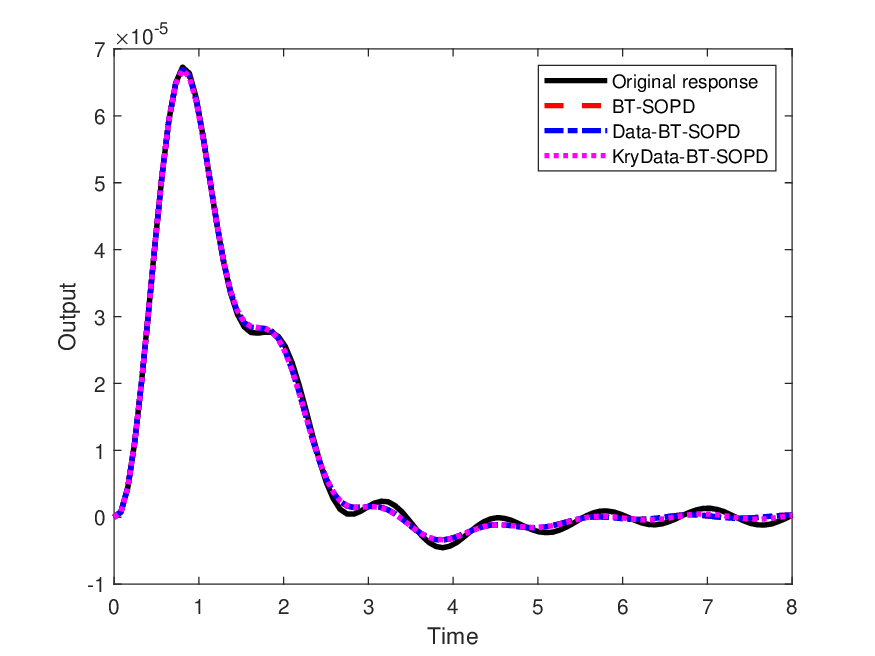}
	\end{minipage}
	\begin{minipage}{0.4\textwidth} 
		\centering 
		\includegraphics[width=0.9\textwidth]{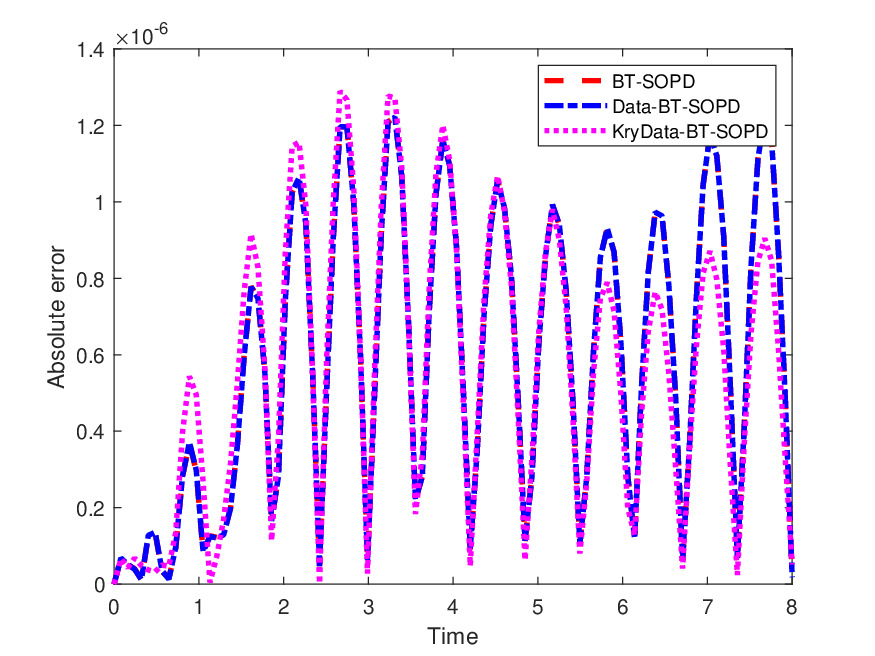}
	\end{minipage}
	\caption{The time response (left) and absolute errors (right) with $r=4, N_q, 
	N_p=500, m=20$ for $u(t)=\mathrm{e}^{-t}\sin(t)$.}
	\label{fig13}
\end{figure}

\begin{table}[htb]
	\centering
	\caption{\centering The CPU time to obtain reduced models with respect to 
	different 
	parameters. }
	\label{tab:11}
	\begin{tabular}{ccccc}
		\toprule
		&$N_q, N_p=50$&$N_q, N_p=500$&$N_q, N_p=500$&$N_q, N_p=500$\\ 
		&$m=5$&$m=5$&$m=20$&$m=30$\\ 
		\midrule
		Data-BT-SOPD & 0.0202s & 0.7460s & 0.7460s & 0.7460s \\
		KryData-BT-SOPD& 0.0281s & 0.0604s & 0.0956s & 0.1366s  \\
		\bottomrule
	\end{tabular}
\end{table}
\end{example}

\begin{example}

This example is a variant of the clamped beam model provided in 
\cite{Benner2005book}, 
which is obtained by spacial discretization of an appropriate partial differential 
equation. The input represents the
force applied to the structure at the free end, and the output is the resulting
displacement. The dimension of the second order system is $174$, and we use 
$\alpha=\beta=0.06$ for the damping matrix $D=\alpha M+\beta K$.

\begin{figure}[htb]
	\centering
	\begin{minipage}{0.4\textwidth} 
		\centering 
		\includegraphics[width=0.9\textwidth]{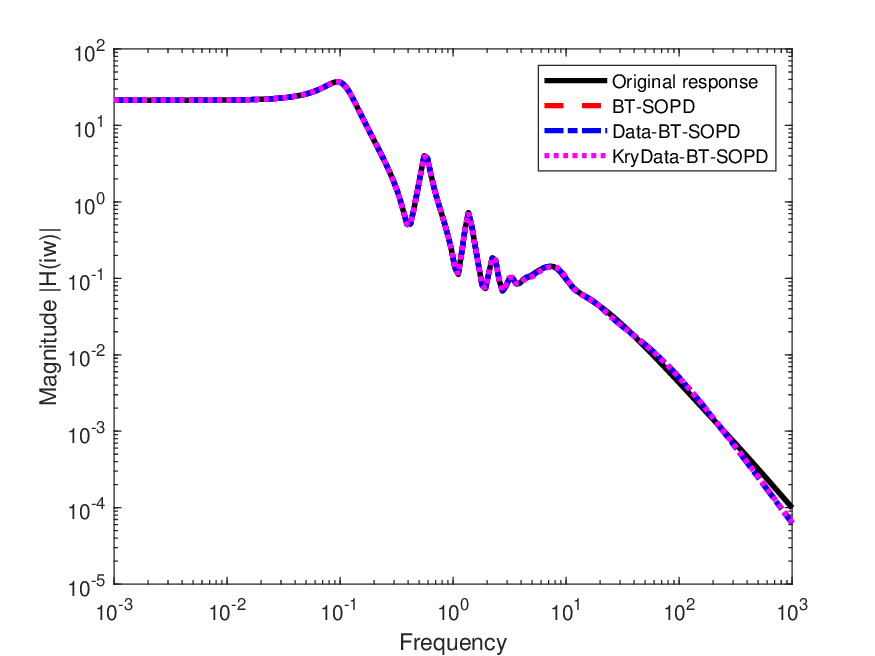}
	\end{minipage}
	\begin{minipage}{0.4\textwidth} 
		\centering 
		\includegraphics[width=0.9\textwidth]{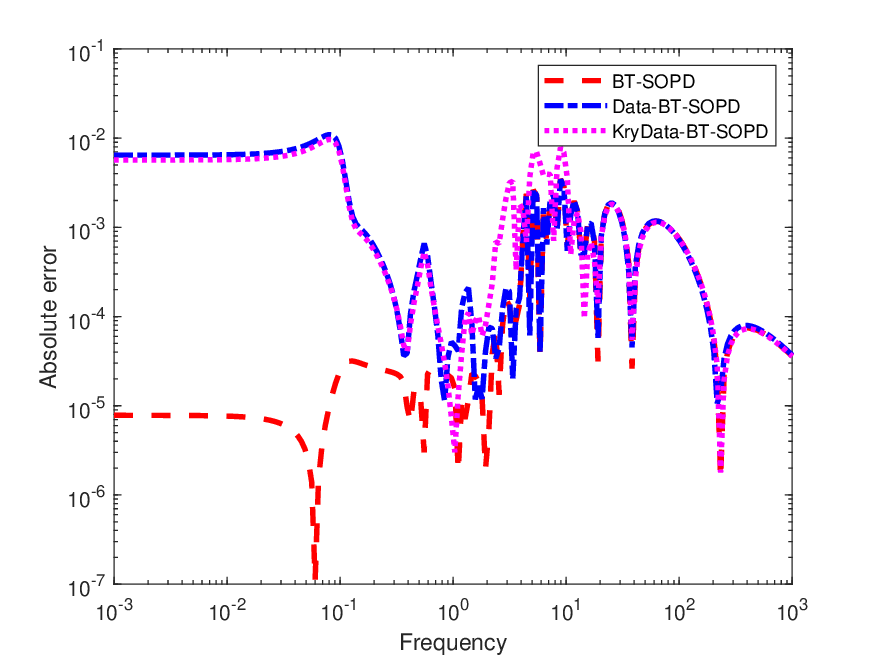}
	\end{minipage}
	\caption{The frequency response (left) and absolute errors (right) with $r=10, 
		N_q, 
		N_p=200, m=20$.}
	\label{fig21}
\end{figure}	
\begin{figure}[htb]
	\centering
	\begin{minipage}{0.4\textwidth} 
		\centering 
		\includegraphics[width=0.9\textwidth]{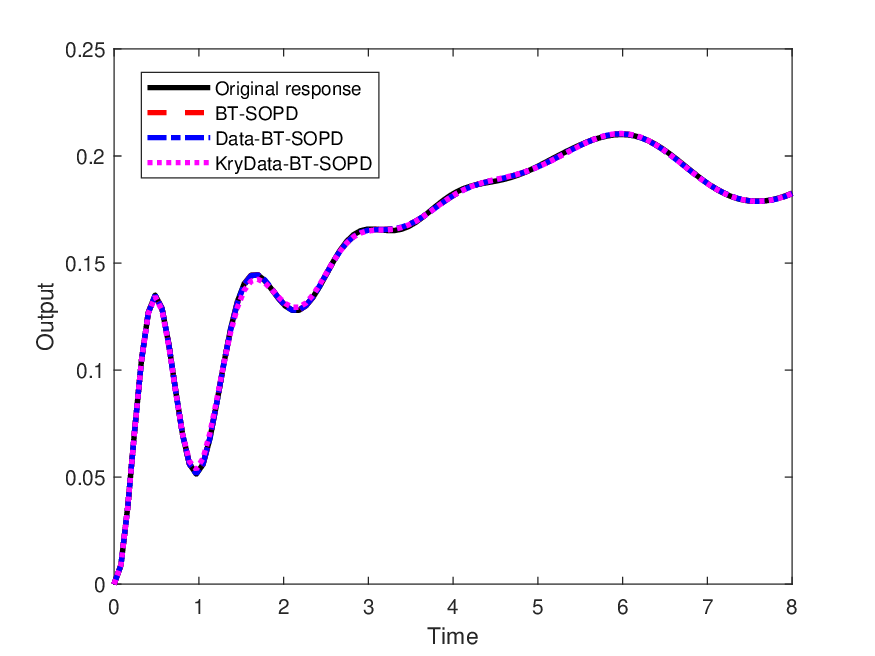}
	\end{minipage}
	\begin{minipage}{0.4\textwidth} 
		\centering 
		\includegraphics[width=0.9\textwidth]{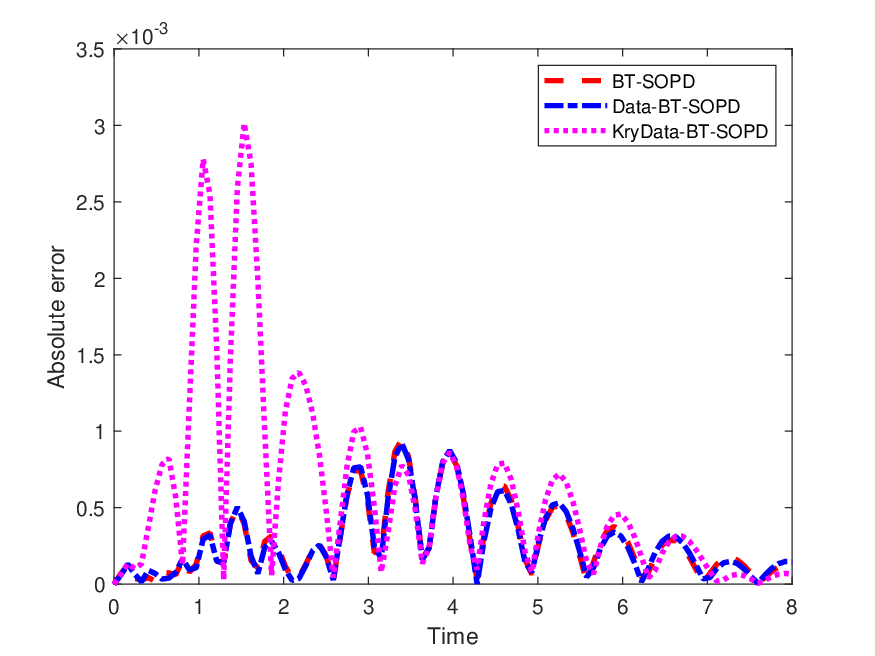}
	\end{minipage}
	\caption{The time response (left) and absolute errors (right) with $r=10, 
		N_q, 
		N_p=200, m=20$ for $\mathrm{e}^{-t}\sin(5t)$.}
	\label{fig22}
\end{figure}

In the simulation, we sample the interval $s\in[10^{-1}, 10^4]$	logarithmically to 
get 200 nodes $s_i$. The sample data is obtained by evaluating the transfer function 
at the points $\mathrm{i} s_i$ and $-\mathrm{i} s_i$. With the reduced order $r=10$, 
Algorithms \ref{SOBT}-\ref{alg:dbt-k} are executed to generate reduced models, where 
$m=20$ is adopted in Algorithm \ref{alg:dbt-k}. Figure \ref{fig21} presents the 
response and the associated absolute errors in the frequency domain. The transient 
behavior of the original system is approximated faithfully by all reduced models, 
and the performance of standard BT is accurately captured by the proposed 
data-driven approach in the high frequency domain for this example. Reduced models 
Data-BT-SOPD and KryData-BT SOPD exhibit a similar performance in this settings. 
Figure \ref{fig22} provides the time response and the according absolute errors for 
the given input $\mathrm{e}^{-t}\sin(5t)$, where Data-BT-SOPD and BT-SOPD almost 
have the same performance in the time domain. While a relatively large error is 
observed for KryData-BT-SOPD in Figure \ref{fig22}, one can 
execute Algorithm \ref{alg:dbt-k} with respect to a higher order extended 
Krylov subspace, say $m=30$, the error of KryData-BT-SOPD in the time domain 
decays notably.

Furthermore, we vary the reduced order from $r=5$ to $r=25$ and test 
the error of each reduced model generated by 
Algorithm \ref{SOBT}-\ref{alg:dbt-k} with the choice $m=30$. Table \ref{tab:21} 
records the relative $H_2$ and $H_\infty$ error of different methods, illustrating 
that Data-BT-SOPD model nearly replicates the performance of BT-SOPD. Generally, the 
order of the extended Krylov subspace in Algorithm \ref{alg:dbt-k} should grow 
as the reduced order and the size 
of samples increase. Because we use a fixed value of $m$ in the simulation, the 
approximation accuracy of KryData-BT-SOPD declines as the reduced order rises in 
Table \ref{tab:21}. For $r=10, m=30$, we also present the CPU time spent on 
constructing reduced models with different sizes of measurements in Table 
\ref{tab:22}. It is clear that KryData-BT-SOPD model is more efficient when a 
amount of samples are involved in the modeling.

\begin{table}[htb]
	\centering
	\caption{\centering Relative $H_2$ and $H_\infty$ errors of reduced models 
	generated by different methods. }
	\label{tab:21}
	\begin{tabular}{cccccc}
		\toprule
		methods&$r=5$&$r=10$&$r=15$&$r=20$&$r=25$\\ 
		\midrule
		BT-SOPD ($H_2$)&2.8674e-02&1.5421e-03&5.5912e-05 &2.5436e-06 &5.2645e-07\\ 
		Data-BT-SOPD 
		($H_2$)&2.8698e-02&1.5599e-03&7.1345e-05&3.9300e-06&3.5216e-07\\ 
		KryData-BT-SOPD ($H_2$)&2.8698e-02&1.5594e-03&1.0533e-04 &7.2890e-05&Inf\\ 
		BT-SOPD ($H_\infty$)&2.9544e-03&9.8511e-05&2.4410e-06 
		&1.9136e-07&3.3826e-08\\
		Data-BT-SOPD ($H_\infty$)&3.7713e-03&3.5141e-04&1.3427e-05&1.5204e-07 
		&1.2761e-08\\ 
		KryData-BT-SOPD 
		($H_\infty$ )&3.7720e-03&3.5164e-04&1.4625e-05&1.2723e-05&1.3677e-05\\ 
		\bottomrule
	\end{tabular}
\end{table}

\begin{table}[htb]
	\centering
	\caption{\centering The CPU time to obtain reduced models with respect to 
	the different choices of $N_q=N_p$. }
	\label{tab:22}
	\begin{tabular}{cccccc}
		\toprule
		methods&$200$&$400$&$600$&$800$&$1000$\\
		\midrule
		Data-BT-SOPD &0.1332s&0.4891s&1.0278s&1.8927s&2.9736s\\ 
		KryData-BT-SOPD &0.0800s&0.1171s&0.1935s &0.3261s&0.4390s\\ 
		\bottomrule
	\end{tabular}
\end{table}
\end{example}

\section{Conclusions} \label{sec:sec-5}

We have investigated the nonintrusive version of BT, namely data-driven BT, for 
second-order systems based on the measurement in the frequency domain. The derived 
relationship between the main quantities in BT and the measurements provides the 
insights on the execution of BT in a nonintrusive manner. The structure-preserving 
reduced models are generated directly based on the frequency domain samples for 
second-order systems with proportional damping, and the execution of the proposed 
approach in real arithmetic is also provided in detail. The low-rank 
approximation to the solution of Sylvester equations enables 
an efficient execution of the proposed approach when a amount of sample data are 
involved in the modeling. The numerical simulation results show that the proposed 
methods nearly replicates the performance of the standard BT methods.


\appendix
\setcounter{figure}{0}
\addcontentsline{toc}{section}{Reference}
\markboth{Reference}{}
\bibliographystyle{elsarticle-num}
\bibliography{reference}

\end{document}